%% file: samp.tex
\documentclass[11pt,letterpaper]{amsart}
\usepackage[utf8]{inputenc}
\usepackage{subcaption}
\usepackage[font=small,labelfont=bf]{caption}
\usepackage{verbatim} 
\usepackage{epsfig,color}
\usepackage{graphicx}
\usepackage{amsmath}
\usepackage{amsfonts}
\usepackage{amssymb}
\usepackage{hyperref}
\newcommand{\R}{\mathbb{R}}
\newcommand{\supp}{\textnormal{supp }}
\newcommand{\WF}{\textnormal{WF}}

\newtheorem{theorem}{Theorem}[section]

\theoremstyle{definition}
\newtheorem{definition}[theorem]{Definition}
\theoremstyle{definition}

\usepackage[backend=bibtex,style=numeric,firstinits=true]{biblatex}
\title[Sampling in TAT]{Sampling in Thermoacoustic Tomography}
\author[C. Mathison]{Chase Mathison}
\date{\today}
\begin{document}
\begin{abstract}
  We explore the effect of sampling rates when measuring data given by
  $Mf$ for special operators $M$ arising in Thermoacoustic Tomography.
  We start with sampling requirements on $M f$ given $f$ satisfying
  certain conditions.  After this we discuss the resolution limit on
  $f$ posed by the sampling rate of $M f$ without assuming any
  conditions on these sampling rates.  Next we discuss aliasing
  artifacts when $M f$ is known to be under sampled in one or more of
  its variables. Finally, we discuss averaging of measurement data and
  resulting aliasing and artifacts, along with a scheme for
  anti-aliasing.
\end{abstract}
\maketitle
\section{Introduction}
This work builds on the theory laid out in \cite{StefanovP2018} on sampling
Fourier Integral Operators (FIOs).  We discuss the specific application of
Thermoacoustic Tomography, in which case the measurement operator $M$ is an FIO
under suitable conditions. We discuss the theoretical resolution of $f$ given
the sampling rate of $M f$ and then discuss aliasing and averaged data. Lastly
we will show empirical evidence of our findings using numerical simulations.

Thermoacoustic Tomography is a medical imaging method in which a short
pulse of electromagnetic radiation is used to excite cells in some
object we wish to image, typically the organs of a patient.  Upon
absorbing the EM radiation, the cells in the patient in turn vibrate,
creating ultrasonic waves that then propagate out of the patient and
are measured by any number of methods.  Using this measured data, we
then try to reconstruct, in some sense, an image of the inside of the
patient.  This is a hybrid imaging method which uses high contrast,
low resolution EM radiation to excite the cells; and low contrast,
high resolution ultrasound waves as measurement
\cite{Oraevsky1994,Kuchment2014,Kruger1999,Kruger2000,wang2015}.  The hope is
to be able to get an image with good contrast and resolution by
combining these two types of waves.

More precisely, let $\Omega \subset \R^n$ be an open subset
of Euclidean $n$-space such that $\bar{\Omega} \subset B_R(0)$ for
some $R > 0$ where $B_R(0)$ is the Euclidean ball of radius $R$.
Suppose $f$ is a smooth function on $\R^n$ supported in $\Omega$.  We
view $f$ as the initial pressure distribution internal to some object
to be imaged. Then, after exposing $\Omega$ to EM radiation, the
ultrasonic waves created solve the acoustic wave equation:
\begin{equation}\label{eq:wave}
\begin{cases}
\left(\partial_t^2 - c^2(x) \Delta_{g_0}\right) u = 0 & (t,x) \in
[0,\infty) \times \R^n, \\
u\mid_{t=0} = f(x) & x \in \R^n, \\
\partial_t u\mid_{t=0} = 0 & x \in \R^n.
\end{cases}
\end{equation}
Here, $c(x) > 0$ is the wave speed, which we take to be identically
$1$ outside of $K \subset \subset \Omega$.  We assume that $c$ is a
smooth function of $x$. In addition, $g_0$ is the Riemannian metric on
the space $\bar{\Omega}$, assumed to be Euclidean on
$\partial \Omega$.  We define $g:= c^{-2}g_0$, which is the metric
form which determines the geometry of this problem.  Assume $u(t,x)$
is a solution to (\ref{eq:wave}) for all
$(t,x) \in [0,\infty) \times \R^n$.  Further suppose that we have
access to $u(t,y)$ for $(t,y) \in (0,T)\times \Gamma$ where $T > 0$
and $\Gamma \subset \partial \Omega$ is a relatively open subset of
$\partial \Omega$ (for this paper, we will take
$\Gamma = \partial\Omega$).  We define for
$(t,y) \in (0,T) \times \Gamma$ the distribution $M f$ as the
measurement operator:
\begin{gather*}
  M: C_0^\infty(\Omega) \rightarrow C_{(0)}^\infty((0,T)\times \Gamma), \\
M f (t,y) = u(t,y),\, (t,y) \in (0,T) \times \Gamma,
\end{gather*}
where $C_{(0)}^\infty((0,T)\times \Gamma)$ is the space of smooth
functions $\phi$ on $(0,T)\times\Gamma$ such that $\phi(t,y) = 0$ near
$t = 0$.  The methods used to collect data on $\Gamma$ are varied and
include point detectors \cite{Kunyansky2008,StefanovP2009,hristova2008,hristova2009}, integrating line
detectors \cite{Burgholzer2006,Grun2007}, circular integrating
detectors \cite{Haltmeier2007,Zangerl2009}, and 2D planar detectors
\cite{Stefanov2017,Haltmeier2004}. We note that at least when
$f \in C_0^\infty(\Omega)$, by energy estimates, $M$ is well defined.
We may actually even take $f$ to be a distribution in
$\mathcal{D}'(\Omega)$ such that
$\lVert f \rVert_{H_D} = \int_{\Omega} |\nabla f|^2\,dx < \infty$, and
by conservation of energy, $M$ extends to a well defined operator.
The closure of $C_0^\infty(\Omega)$ under the previously stated norm
is the space $H_D(\Omega)\subset H_0^1(\Omega)$, and we will assume $f \in H_D(\Omega)$
unless otherwise stated.

\subsection{\(M\) as an FIO}

To obtain an oscillatory integral representation of $M$, we may use
the geometric optics construction to solve for $u(t,x)$ in
$(0,T)\times \R^n$ up to a smooth error (see
\cite{StefanovP2009,Taylor81} for more details).  This construction
leads to the representation
\begin{gather*}
  u(t,x) = \frac{1}{(2\pi)^{n}} \sum_{\sigma = \pm} \int
  e^{i\phi_{\sigma}(t,x,\xi)}a_{\sigma}(t,x,\xi)\hat{f}(\xi)\,d\xi,
\end{gather*}
where $\phi_\sigma$ are solutions to the eikonal equation
$\left(\partial_t \phi_\sigma\right)^2 =c^2(x)|\nabla_x \phi_\sigma
|_{g_0}^2$ with initial conditions
$\phi_\sigma (0,x,\xi) = x\cdot \xi$.  Note that solutions to the
eikonal equation are local in nature, and so this representation of 
$u(t,x)$ is only valid until some time $T_1$.  However, we may then
solve (\ref{eq:wave}) with ``initial'' conditions
$\tilde{u}(0,x) = u(T_1,x)$ and
$\partial_t \tilde{u}(t,x)\mid_{t=0} = \partial_t u(t,x)\mid_{t=T_1}$
using the same geometric optics construction.  In this way, we can
obtain an ``approximate'' solution to (\ref{eq:wave}) for all $(t,x)$.
Note by approximate, we mean up to a smooth error term.  This error
term could be quite large in the $L^\infty$ sense, but because it is a
smooth term, it is negligible in the calculus of FIOs. It
can be shown that $M = M_+ + M_-$ is a sum of elliptic FIOs of order
$0$ associated with locally diffeomorphic canonical relations that are
each (locally) one-to-one mappings (see
i.e. \cite{StefanovP2018,StefanovP2009}).  We record the canonical
relations $C_+$ and $C_-$ here for later use:
\begin{gather}\label{eq:canonical_rel}
  C_{\pm}: \left(x,\xi\right) \mapsto
  \left(s_{\pm}(x,\xi),\gamma_{x,\xi}(s_{\pm}(x,\xi)),\mp
    |\xi|_g,\dot{\gamma}_{x,\xi}'(s_{\pm}(x,\xi))\right).
\end{gather}
Here, we have $s_{\pm}(x,\xi)$ is the exit time of the geodesic
starting at $x$ in the direction $\pm g^{-1}\xi$, $\gamma_{x,\xi}(t)$
is the point on the geodesic issued from $(x,\xi)$ at time $t$ and
$\dot{\gamma}_{x,\xi}'(t)$ is the orthogonal (in the metric)
projection of $\dot{\gamma}_{x,\xi}(t)$ onto $T\partial\Omega$ (the
tangent bundle of the boundary of $\Omega$, so implicitly, we assume
that $\partial \Omega$ is a at least a $C^1$ manifold).  We assume
that the metric induced by $g:=c^{-2}(x)g_0$ is non trapping, so that
$|s_{\pm}(x,\xi)| < \infty$ for all $(x,\xi) \in T^*\Omega$.  Note
that because each of the canonical relations $C_+$ and $C_-$ are
one-to-one, the full canonical relation of the FIO $M$ given by
$C = C_+ \cup C_-$ is one-to-two, which makes intuitive sense as
singularities split and travel along geodesics according to
propagation of singularities theory.
\section*{Acknowledgments}
The author would like to thank Dr. Plamen Stefanov for suggesting this
problem and for his guidance in the analysis of this problem.
\section{Preliminary definitions and theorems}
\subsection{Semiclassical analysis}
The main definitions and theorems of semiclassical analysis and
sampling that we use come from \cite{Zworski2012,StefanovP2018}.  For
a more complete background on semiclassical analysis, see
\cite{Zworski2012}. In sampling the measurement operator $M f$, we are
interested in how the sampling rates affect our ability to resolve
singularities with high frequency.  To model this, we will rescale
co-vectors $\xi$ by a factor of $1/h$ where $h$ is a small parameter.
We then examine families of functions (or distributions) $f_h$ that
satisfy certain growth conditions as $h$ becomes small.  Because of
this, instead of considering the classical wave front set of a
distribution, we consider the semiclassical wave front set, denoted
$\WF_h(f)$.  Note that $f$ is understood here to be a family of
functions $f_h$ depending on the parameter $h$, but we will drop this
subscript when it will not cause confusion.  A key tool in analyzing
the behavior of the measurement operator $M$ will be the semiclassical
Fourier Transform, defined below.
\begin{definition}[Semiclassical Fourier Transform]\label{def:scft}
The semiclassical Fourier transform of an $h$-dependent family of distributions is defined as
\begin{gather*}
\mathcal{F}_h f_h(\xi) = \int e^{-ix\cdot \xi/h} f_h(x)\,dx.
\end{gather*}
If we denote the classical Fourier Transform by $\mathcal{F}$, then we have
\begin{gather*}
\mathcal{F}_h f_h(\xi) = \mathcal{F} f_h\left(\frac{\xi}{h}\right).
\end{gather*}
\end{definition}
Much like in classical analysis, we can use the semiclassical Fourier
transform to define Sobolev norms on certain classes of functions or
distributions.
\begin{definition}[$h$-Tempered family of
  distributions]\label{def:htempered}
  The $h$-dependent family $f_h$ of distributions in $\mathcal{S}'$ is
  said to be $h$-tempered if
\begin{gather*}
  \lVert f_h \rVert_{H_h^s}^2 := (2\pi h)^{-n} \int \langle \xi
  \rangle^{2s} \left| \mathcal{F}_h f(\xi) \right|^2\,d\xi
\end{gather*}
is such that $\lVert f_h \rVert_{H_h^s} = O(h^{-N})$ for some $s$ and
$N$.  Here, we have $\langle \xi \rangle = \sqrt{1 + |\xi|^2}$.
\end{definition}

Another key tool we will use is the idea of the semiclassical wave
front set of an $h$-dependent family of distributions.
\begin{definition}[Semiclassical Wave Front Set]\label{def:scwfs}
  The semiclassical wave front set $\WF_h(f_h)$ of the $h$-tempered
  family $f_h$ is defined to be the complement of the set of
  $ (x_0,\xi^0) \in \R^{2n}$ such that there exists
  $\phi \in C_0^\infty(\R^n)$ with $\phi(x_0) \neq 0$ so that
\begin{gather*}
\mathcal{F}_h (\phi f_h) = O(h^\infty)\qquad \text{(in }L^\infty\text{)}
\end{gather*}
for $\xi$ in a neighborhood of $\xi^0$.
\end{definition}
This set plays a similar role as the classical wave front set from
microlocal analysis, however in general there is no sort of inclusion
between these two sets.  As an example \cite{Zworski2012}, the
coherent state
$$f_h(x) = e^{-|x-x_0|^2/(2h)}e^{ix\cdot \xi^0/h},$$
has an empty wave front set in the classical sense, as it is a smooth
function in both $x$ and $\xi$, however its semiclassical wave front
set is $\WF_h(f_h) = \lbrace (x_0,\xi^0)\rbrace$.  Note also that the
zero section is allowed to be a part of the semiclassical wave front
set, unlike in the classical case.  Also, we do not require the
semiclassical wave front set to be a conic set, which is another way
that this set differs from the classical wave front set.

We call elements of $\WF_h(f_h)$ singularities, even though a function
with finite semiclassical wave front set is actually smooth.

\begin{definition}[$h$-$\Psi$DO]\label{def:hpdo}
  We will use the standard quantization to define semiclassical
  pseudodifferential operators.  Fix $m$ and $k\in \R$ and let
  $a(x,\xi) \in C^{\infty}(\R^{2n})$ satisfy the following: For every
  $\alpha$ and $\beta$ multi-indices and every compact set
  $K \subset \R^n$ there exists some $C_{\alpha,\beta,K} > 0$ such
  that
$$
|D_x^{\alpha} D_{\xi}^{\beta} a(x,\xi)| \le C_{\alpha,\beta,K} h^k \langle \xi\rangle^{m}
$$
for all $x\in K$ and $\xi \in \R^n$.  We then say $a(x,\xi)$ is a
semiclassical symbol of order $\le m$.  Then we define the
semiclassical pseudodifferential operator $a(x,hD)$ by
$$
a(x,hD)f(x) := (2\pi h)^{-n} \iint_{\R^{2n}} e^{i(x-y)\cdot \xi /
  h}a(x,\xi) f(y)\,dy\,d\xi.
$$
\end{definition}
\begin{definition}\label{def:phloc}
  The $h$-tempered family $f_h$ is said to be localized in phase space
  if there exists some $\psi \in C_0^{\infty}(\R^{2n})$ such that
\begin{gather*}
\left(\text{Id} - \psi(x,hD)\right)f_h = O_{\mathcal{S}}(h^{\infty}).
\end{gather*}

Note that because the functions we work with are semiclassically band
limited (see definition \ref{def:scbl}), that all functions we work
with can be assumed to be localized in phase space unless otherwise
stated.
\end{definition}

\begin{definition}[Semiclassical Frequency Set]\label{def:scfs}
For each tempered $h$-dependent distribution $f_h$ localized in phase space, set
\begin{gather*}
\Sigma_h(f_h) = \left\{ \xi \mid (x,\xi) \in \WF_h (f_h) \text{ for some } x\in \R^n \right\}.
\end{gather*}
This is simply the projection of $\WF_h(f_h)$ onto the second variable.
\end{definition}

\begin{definition}[Semiclassically Band Limited Functions]\label{def:scbl}
We say that $f_h \in C_0^\infty(\R^n)$ is semiclassically band limited (in $\mathcal{B}$) if\begin{enumerate}
\item $\supp f_h$ is contained in an $h$-independent set,
\item $f_h$ is tempered,
\item there exists a compact set $\mathcal{B}\subset \R^n$ such that
  for every open $U\supset \mathcal{B}$, we have for every $N$ there
  exists $C_N$ such that
\begin{gather*}
|\mathcal{F}_h f_h (\xi)| \le C_N h^N \langle \xi \rangle^{-N}\, \text{for } \xi \not\in U.
\end{gather*}
\end{enumerate}
\end{definition}

Semiclassically band limited functions are those functions that can be
reconstructed up to a smooth error from their samples, much like the
band limited functions are those that can be perfectly reconstructed
from their samples in the classical Nyquist Sampling theorem given a
small enough sampling rate\cite{Marks1991}.
\subsection{Sampling}
The main theorem used in
\cite{StefanovP2018} is the following:
\begin{theorem}\label{thm:samp}
  Assume that $\Omega \subset \R^n$, $\mathcal{B} \subset \R^n$ are
  open and bounded.  Let $f_h \in C_0^\infty(\Omega)$ satisfy
\begin{gather}
  \lVert (Id - \psi(x,hD))f_h \rVert_{H_h^m} = O(h^\infty) \lVert
  f_h\rVert,\qquad \forall m \gg 0,
\end{gather}
for some $\psi \in C_0^\infty(\R^{2n})$ such that
$\text{supp}_\xi \psi \subset \mathcal{B}$.  Let
$\hat{\chi}\in L^\infty(\R^n)$ be such that
$\supp \hat{\chi} \subset \mathcal{B}$ and $\hat{\chi}=1$ near
$\text{supp}_\xi \psi$.

Assume that $W$ is an invertible matrix so that the images of
$\mathcal{B}$ under the translations
$\xi \mapsto \xi + 2\pi(W^*)^{-1}k,\, k\in\mathbb{Z}^n$, are mutually
disjoint.  Then for every $s\in (0,1]$,
\begin{gather}
  f_h(x) = |\det W| \sum_{k\in \mathbb{Z}^n} f_h(shWk)\chi\left(
    \frac{\pi}{sh}(x-shWk)\right) + O_{H^m}(h^\infty)\lVert f_h
  \rVert_{L^2},
\end{gather}
for every $m \ge 0$, and
\begin{gather}
  \lVert f_h \rVert_{L^2}^2 = |\det W|(sh)^n \sum_{k\in \mathbb{Z}^n}
  |f_h(shWk)|^2 + O(h^\infty)\lVert f \rVert_{L^2}^2.
\end{gather}
\end{theorem}

The proof of this theorem essentially follows from the classical
Nyquist sampling theorem and can be found in \cite{StefanovP2018,Petersen1962}.
For all applications in this paper, we take the matrix $W$ above to be
the identity matrix.

We make heavy use of the following theorem which relates how classical
FIOs effect semiclassical wavefront sets from \cite{StefanovP2018},
where the reader can find the proof.
\begin{theorem}\label{thm:canon}
  Let $A$ be an FIO in the class $I^m(\R^{n_2},\R^{n_1},\Lambda)$
  where $\Lambda \subset T^*(\R^{n_1}\times \R^{n_2})\setminus 0$ is a
  Lagrangian manifold and $m \in \R$.  Then for every $f_h$ localized
  in phase space,
\begin{gather}
\WF_h(Af)\setminus 0 \subset C \circ \WF_h(f) \setminus 0,
\end{gather}
where $C = \Lambda'$ is the canonical relation of $A$.
\end{theorem}

This theorem shows how classical FIOs affect the semiclassical
wavefront set away from the zero section. In particular, the
semiclassical wavefront set of $Af$ away from the zero section
transforms in the same way the classical wavefront set does: it is
transformed by the canonical relation associated with $A$.
The main assertion in \cite{StefanovP2018} is that the sampling
requirements of $Mf$ given $\WF (f)$ are determined by $C$, the
canonical relation associated with $Mf$.
\section{Resolution limit of \(f\) given sampling rate of \(M f\)}
Suppose we wish to sample the $M f$ at some fixed sampling rates $s_t$
and $s_{y^j}$.  Here we don't assume that we know any information about
$\Sigma_h(f)$, we only wish to see how fixing a sampling rate on $M f$
affects our ability to resolve singularities of $f$.  Avoiding
aliasing of $M f$ is equivalent to (by Theorem \ref{thm:samp})
\begin{gather*}
(\tau,\eta)\in \Sigma_h(M f)\implies |\tau| \le \frac{\pi}{s_t}, \qquad |\eta_j| \le \frac{\pi}{s_{y^j}},
\end{gather*}
where $\tau$ is the dual variable to $t$, and $\eta$ is the dual
variable to $y$, with $\eta_j$ the $j$th component of $\eta$.  Note that the norms $|\tau|$ and $|\eta|$ are taken
in the corresponding metric. In particular, although the norm $g_0$ on
$\bar{\Omega}$ is assumed to be Euclidean, the induced norm on the
tangent space to the boundary, which we'll call $g_{0,\partial \Omega}$, is not
necessarily Euclidean.  We may use the canonical relation
(\ref{eq:canonical_rel}) $C$ associated with $M$ to write the
inequalities above as
\begin{gather*}
  |\xi|_g = \sqrt{c^2 g_0^{ij}\xi_i\xi_j} \le \frac{\pi}{s_t},\qquad
  |\dot{\gamma}_{x,\xi}'(s_\pm(x,\xi))_j|_{g_{0,\partial\Omega}} \le \frac{\pi}{s_{y^j}}.
\end{gather*}
From this we can see that we have that avoiding aliasing is equivalent to
\begin{gather}\label{eq:resIneq}
c(x) |\xi|_{g_0} \le \frac{\pi}{s_t},\qquad
|\dot{\gamma}_{x,\xi}'(s_{\pm}(x,\xi))_j|_{g_{0,\partial\Omega}} \le \frac{\pi}{s_{y^j}}
\end{gather}
For most of the paper, we will assume that $g_0$ is Euclidean,
although more general results hold.
\subsection{The effect of \texorpdfstring{$s_t$}{s\_ t} on resolution}
Consider the first inequality in (\ref{eq:resIneq}) and assume that
$s_{y^j}$ is taken small enough so as to not effect resolution of
singularities of $f$. The first inequality indicates that the sampling
rate $s_t$ imposes a limit on the resolution of $f$ such that for
fixed $x$, there will be higher resolution of singularities of $f$ at
points $(x,\xi)$ where the wave speed $c(x)$ is slower, and likewise
the resolution will be worse at those points $(x,\xi)$ where the wave
speed is faster. In particular, given the relative sampling rate
$s_t$, we cannot resolve singularities at $x$ with frequency greater
than $$|\xi| = \frac{\pi}{c(x)s_t}.$$ This is a local result. A global
estimate for the maximum frequency of a singularity that is guaranteed
to be resolved anywhere given the sampling rate $s_t$ is given by
\begin{equation}\label{eq:resSt}
|\xi| = \frac{\pi}{c_{\max} s_t}.
\end{equation}
This is illustrated in Figures \ref{fig:res_fast_t_v2} and
\ref{fig:res_slow_t_v2} below.
\begin{figure}[!ht]
\centering
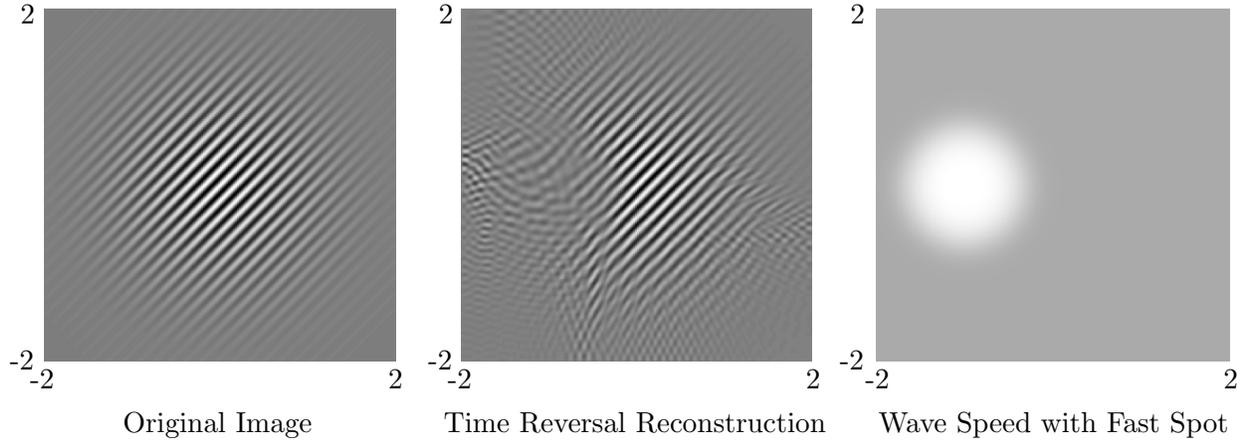
\caption{Resolution of $f$ given a fixed sampling rate $s_t$ of
  $Mf(t,y)$.  The wave speed here
  $c(x,y) = 1 + 0.5\exp(-((x+1)^2 + y^2)^2/0.25)$ has a fast spot
  centered at $x=-1$.  We can see that this is precisely where the
  reconstruction of $f$ has poor resolution when under sampled in the
  $t$ variable, as explained above.}
\label{fig:res_fast_t_v2}
\end{figure}
\begin{figure}[!ht]
\centering
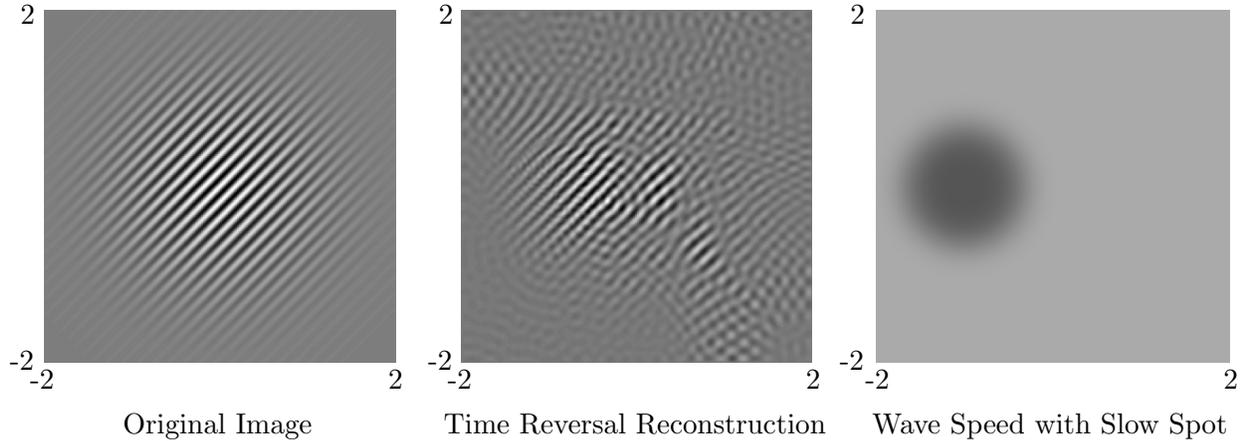
\caption{Resolution of $f$ given a fixed sampling rate $s_t$ of
  $Mf(t,y)$.  The wave speed here
  $c(x,y) = 1 - 0.5\exp(-((x+1)^2 + y^2)^2/0.25)$ has a slow spot
  centered at $x=-1$.  We can see that this is precisely where the
  reconstruction of $f$ has the best resolution when under sampled in
  the $t$ variable, as explained above.}
\label{fig:res_slow_t_v2}
\end{figure}
\subsection{The effect of \texorpdfstring{$s_{y^j}$}{s\_ yj} on resolution}
Assume now that $s_t$ is chosen small enough so as to not effect
resolution of singularities of $f$.  The second inequality in
(\ref{eq:resIneq})
$$
|\dot{\gamma}'_{x,\xi}(s_{\pm}(x,\xi))_j|_{g_{0,\partial\Omega}}\le \frac{\pi}{s_{y^j}},
$$
tells us that the sampling rate $s_{y^j}$ imposes a limit on the
resolution of $f$ such that singularities $(x,\xi)$ that intersect the
boundary $\partial \Omega$ nearly perpendicularly will have higher
resolution than those that hit the boundary nearly tangentially (at a
large angle to the normal vector to $\partial \Omega$ at the point of
intersection).  Also, because $|\dot{\gamma}_{x,\xi}(t)|_g$ is
constant along the geodesic $\gamma_{x,\xi}$, we know in particular
that
$|\dot{\gamma}'_{x,\xi}(s_{\pm}(x,\xi))_j|_{g_{0,\partial\Omega}}\le  |\dot{\gamma}'_{x,\xi}(s_{\pm}(x,\xi))|_{g_{0,\partial\Omega}} = |\xi|_g\cos(\theta)$
where $\theta$ is the angle (in the metric) between
$\dot{\gamma}_{x,\xi}(s_{\pm}(x,\xi))$ and
$\dot{\gamma}'_{x,\xi}(s_{\pm}(x,\xi))$.  This tells us that to avoid
aliasing, we must have
$$
|\xi|_g\cos(\theta) \le \frac{\pi}{s_{y^j}}.
$$
We recall that $|\xi|^2_g = c^2(x)g_0^{ij}\xi_i\xi_j$, and in the case
that $g_0$ is Euclidean, we get
$$
c(x)|\xi| \cos(\theta) \le \frac{\pi}{s_{y^j}}.
$$
For a fixed relative sampling rate $s_{y^j}$, we cannot resolve
singularities $(x,\xi)$ of $f$ of frequency greater than
$$
|\xi| = \frac{\pi}{s_{y^j} c(x) \cos(\theta)}.
$$
Note in particular that if $\theta = \frac{\pi}{2}$ (i.e. the geodesic
$\gamma_{x,\xi}$ hits the boundary $\partial \Omega$ perpendicularly),
then $c(x)|\xi|\cos(\theta) = 0 < \pi/s_{y^j}$, and we will always be able
to resolve the singularity at $(x,\xi)$.  Also note that this is a
local result, and as is the case for $s_t$ ``slow spots'' in the speed
$c(x)$ give better resolution of singularities in general.  Because
$c(x) \le c_{\max}$, we also get the following estimate for the
maximum frequency of a resolvable singularity, regardless of location:
$$
|\xi| = \frac{\pi}{c_{\max}s_{y^j} \cos(\theta)}.
$$
Finally, because $0 < \theta \le \pi/2$, we know
$0 \le \cos(\theta) < 1$, and we have the following (worst case)
global estimate for the maximum frequency of a singularity of $f$ that
can be resolved:
\begin{equation}\label{eq:resSy}
|\xi| = \frac{\pi}{c_{\max}s_{y^j}}.
\end{equation}

\begin{figure}[!ht]
\centering
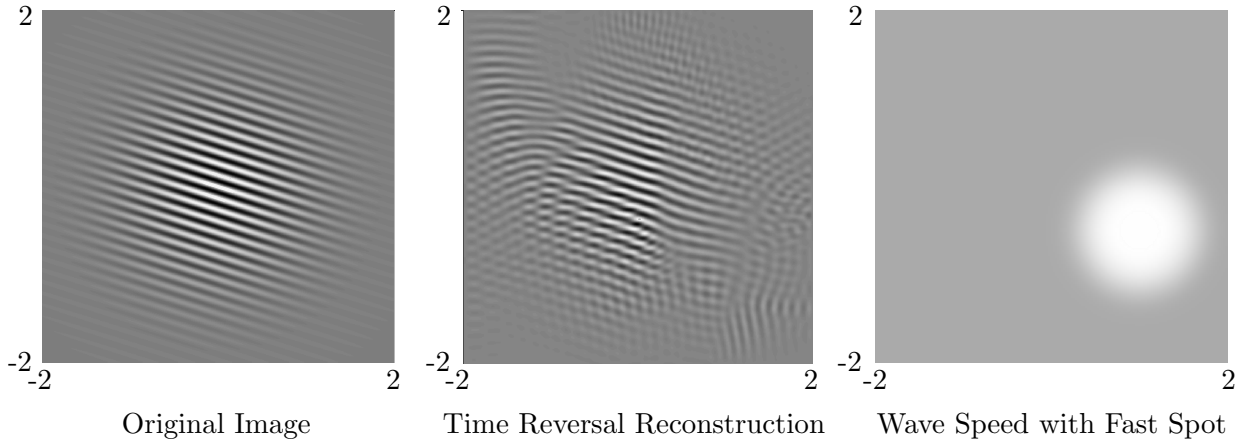
\caption{Resolution of $f$ given a fixed sampling rate $s_{y^j}$ of the
  space variables on the boundary $\partial \Omega$.  We can see that
  the blurring effect is roughly uniform for points near the fast spot
  in the wave speed
  $c(x,y) = 1 + 0.5\exp(-((x-1)^2+(y+0.5)^2)^2)/0.25)$, but that there
  are singularities in the region where $c \approx 1$ far from the
  fast spot that are also highly affected.  These singularities hit the
  boundary with a larger angle to the outward pointing normal vector,
  and so we expect lower resolution there.}
\label{fig:res_fast_y}
\end{figure}

We note that if one wants to be able to resolve singularities of $f$
with frequency $K$, then by considering (\ref{eq:resSt}) and
(\ref{eq:resSy}), the sampling rates $s_t$ and $s_{y^j}$ of $M f$ should
be taken to be at least
$$
s_t = s_{y^j} \le \frac{\pi}{Kc_{\max}},
$$
where $c_{\max}$ is defined as before.  In particular, we recover the
result from \cite{StefanovP2018} that for a semiclassically band
limited $f_h$ with essential maximum frequency $B$ in the Euclidean
case that we need to take sampling rates of $M f$ satisfying
$$
s_t \le \frac{\pi}{B c_{\max}},\qquad s_{y^j}\le \frac{\pi}{B c_{\max}},
$$
to avoid aliasing.  These effects are shown in Figure \ref{fig:res_fast_y}.

\subsection{CFL condition}
We can relate this analysis to numerical solvers of the wave equation.
When solving the wave equation numerically, a typical approach is to
discretize the space and time domain, and use a finite difference
scheme.  Suppose we wish to simulate an experiment using a rectangular
grid in the space coordinates and we collect data on the boundary of a
square. Further, we assume that $g_0$ is Euclidean, and because the
boundary is a rectangle, also the metric induced on the boundary is
Euclidean.  Suppose we have fixed each
$s_{x^j} = \Delta x^j/h \le \pi/(Bc_{\max})$ with a common value
$s_x = \Delta x / h$, where $B$ is the essential band limit on $f$,
i.e. $\Sigma_h(f) \subset [-B,B]^n$.  Note that by our choice of
$s_x$, there will not be aliasing of $Mf$, provided $s_t$ is chosen
well, as on the boundary in this rectangular grid, we have
$s_y = s_x$, where all of the $s_{y^j}$ as above have a common fixed
step size $s_y$. In order to choose $s_t$, we recall that the
frequency set $\Sigma_h (Mf)$ is contained in the set
$\left\{ (\tau,\eta) \mid |\eta| \leq |\tau| \right\}$.  Because $f$
has a semiclassical band limit of $B$, we know that
$\pi_2(\Sigma_h(Mf)) \subset \{|\eta| \leq \sqrt{n} Bc_{\max}\}$,
where $\pi_2$ is the projection onto the second factor. We know this
because each $|\eta_j| \leq Bc_{\max}$. Also, by the analysis above,
we know that $|\tau| = |\xi|_g$, but
$|\xi|_g \le \max |\xi| c_{\max}$.  We also know that
$\max |\xi| < \sqrt{n}B$, so that the largest possible size of
$|\tau|$ given the band limit on $f$, is $\sqrt{n}Bc_{\max}$.  It is
then clear that we need $s_t \leq \pi /(\sqrt{n}B c_{\max}))$ to avoid
aliasing.  This tells us that we should take
$\Delta t \leq \pi h/(\sqrt{n}B c_{\max}) = \Delta x/\sqrt{n}$.  Now,
the CFL condition for the leapfrog finite difference scheme
(\cite{CFL59,Bartels2016,StrangNotes}) tells us that given a step size
$\Delta x$ and wave speed $c(x)$, that we should take the time step
$\Delta t \le \Delta x / (\sqrt{n} c_{\max})$ to ensure stability of
the finite difference scheme.  But
$\Delta x / (\sqrt{n} c_{\max}) \le \Delta x/\sqrt{n}$, because
$c_{\max} \ge 1$.  This means, that if we've chosen
$\Delta x \le \pi h/(B c_{\max})$, and we choose $\Delta t$ satisfying
the CFL condition for the leapfrog finite difference scheme, then
there will be no aliasing in the measured data $Mf$ at the
boundary. Also, if $c_{\max} = 1$, then the CFL condition is identical
to the conditions on $\Delta x$ and $\Delta t$ required to avoid
aliasing of the measured data $Mf$.

\section{Aliasing and artifacts}
Now suppose that we know that $f_h$ is a semiclassically band limited
function with essential band limit $B$.  In \cite{StefanovP2018}, it
is shown that in order to avoid aliasing of $M f_h$, for a
semiclassically band limited $f_h$, we must have relative sample rates
of $s_t \le \frac{\pi N}{B}$ and $s_{y^j} \le \frac{\pi N N'}{B}$
where $B$ is half the side length of a box bounding $\Sigma_h (f)$,
$N$ is the sharp lower bound of the metric form $g=c^{-2}g_0$ on the
unit sphere for all $x$, and $(N')^2$ is the sharp upper bound on the
induced metric on the Euclidean sphere in a fixed chart for $y$. In
the numerical examples that follow, $\partial \Omega$ is piecewise
flat and parameterized in a Euclidean way, so that $N' = 1$ away from
corners. Note that if $g_0$ is Euclidean, then setting
$c_{\max} =\max c(x)$, we have $N = 1/c_{\max}$, and $N' = 1$ so that
the relative sampling rates needed to avoid aliasing are
$$
s_t \le \frac{\pi}{B c_{\max}},\qquad s_{y^j} \le \frac{\pi}{B c_{\max}}.
$$

\subsection{Under sampling in \(t\)}
\begin{figure}[!ht]
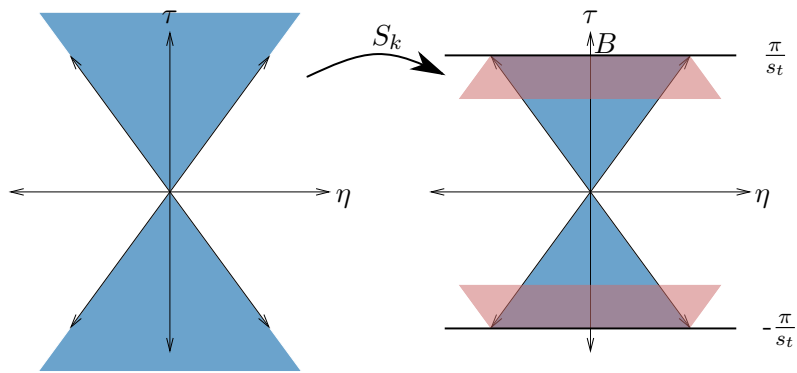
\caption{The characteristic cone in which $\Sigma_h(Mf)$ must lie.
  The cone on the left shows the possible range of the covector
  $(\eta,\tau)$ which is determined by the canonical relation
  associated with $M$.  The image on the right shows the possible
  range of covectors $(\eta,\tau)$ after under sampling (in $t$).  Note
  that the red regions have been shifted up and down from the original
  frequency set by translation due to
  under sampling.}\label{fig:charCone}
\end{figure}

Suppose that we have chosen $s_t$ such that
$s_t > \frac{\pi}{Bc_{\max}}$.  Then, by \cite{StefanovP2018} there
will be aliasing of $M f$.  The error in the reconstruction can be
modeled by the frequency shift operator
$$S_k: \tau \rightarrow \tau + \frac{2\pi k}{s_t}.$$ This operator is
valid as long as $\tau + 2\pi k/s_t \in [-\pi/s_t,\pi/s_t]$ (see
Figure \ref{fig:charCone} (right)).
If we have not under sampled $M f$ too critically in the $t$ variable,
we would expect to only see this added error for $k = -1,1$, with more
terms added as the under sampling becomes worse. As explained in
\cite{StefanovP2018}, by Egorov's Theorem, we expect to see artifacts
in a reconstruction of $f$ that can be calculated by the canonical
relation
\begin{gather*}
C_{\pm}^{-1} \circ S_k \circ C_{\pm}: (x,\xi) \rightarrow (\tilde{x},\tilde{\xi}),
\end{gather*}
where $\tilde{x}$ and $\tilde{\xi}$ can be calculated by finding the
operator on the left.  We do that now for $C_+$:
\begin{align*}
  C_+^{-1}\circ S_k \circ C_+(x,\xi) &= C_+^{-1}\circ S_k (s_+(x,\xi),\gamma_{x,\xi}(s_+(x,\xi)),-|\xi|_g,\dot{\gamma}_{x,\xi}'(s_+(x,\xi))) \\
                                     &= C_+^{-1} (s_+(x,\xi),\gamma_{x,\xi}(s_+(x,\xi)),-|\xi|_g + \frac{2\pi k}{s_t},\dot{\gamma}_{x,\xi}'(s_+(x,\xi))) \\
                                     &= \left(\gamma_{y,-\zeta}(s_+(x,\xi)), -\dot{\gamma}_{y,-\zeta}(s_+(x,\xi))\right),
\end{align*}
where $y = \gamma_{x,\xi}(s_+(x,\xi))$ is the point of intersection of
the geodesic issued from $(x,\xi)$ with $\partial \Omega$, and
$\zeta = \dot{\gamma}'_{x,\xi}(s_+(x,\xi)) + \beta_k \eta^\bot$ where
$\beta_k = \sqrt{ (|\xi|_g - 2\pi k/s_t)^2 -
  |\dot{\gamma}'_{x,\xi}(s_+(x,\xi))|^2}$ and
$\eta^\bot = \dot{\gamma}_{x,\xi}(s_+(x,\xi)) -
\dot{\gamma}'_{x,\xi}(s_+(x,\xi))$.  Aliasing artifacts are found
using this mapping in Figures \ref{fig:sing_const_t} and
\ref{fig:sing_var_t} below.  The mapping $C_-^{-1}\circ S_k \circ C_-$
is calculated in almost an identical fashion, however we have a change
in sign in the $\tau$ variable.
\begin{figure}[!ht]
\centering
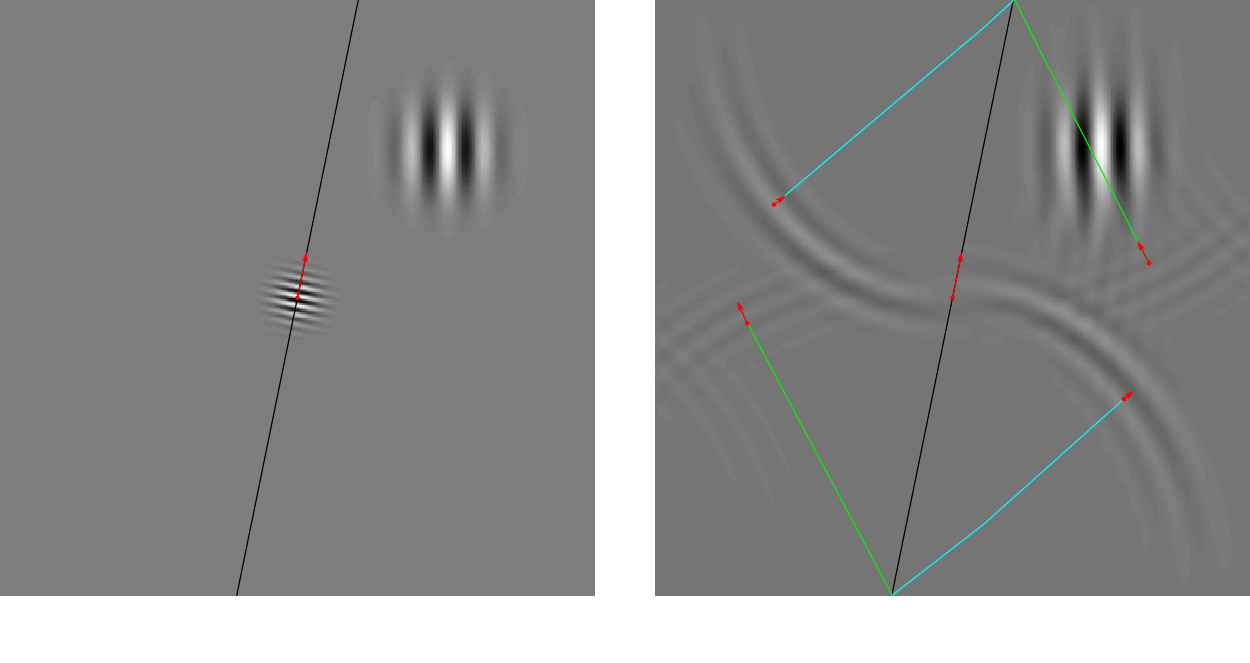
\caption{Tracing the aliasing artifacts by using geodesics. We have
  used the constant wave speed $c \equiv 1$ for this example.  Here we
  have under sampled in $t$ and show the image of the singularity
  $(x,\xi)$ under the canonical relations given by
  $C_{\pm}^{-1} \circ S_i \circ C_{\pm}$ for $i = 1,2$. Note that the
  low frequency singularity does not cause artifact, but the high
  frequency singularity vanishes in the reconstruction and causes
  aliasing artifacts.}
\label{fig:sing_const_t}
\end{figure}

\begin{figure}[!ht]
\begin{center}
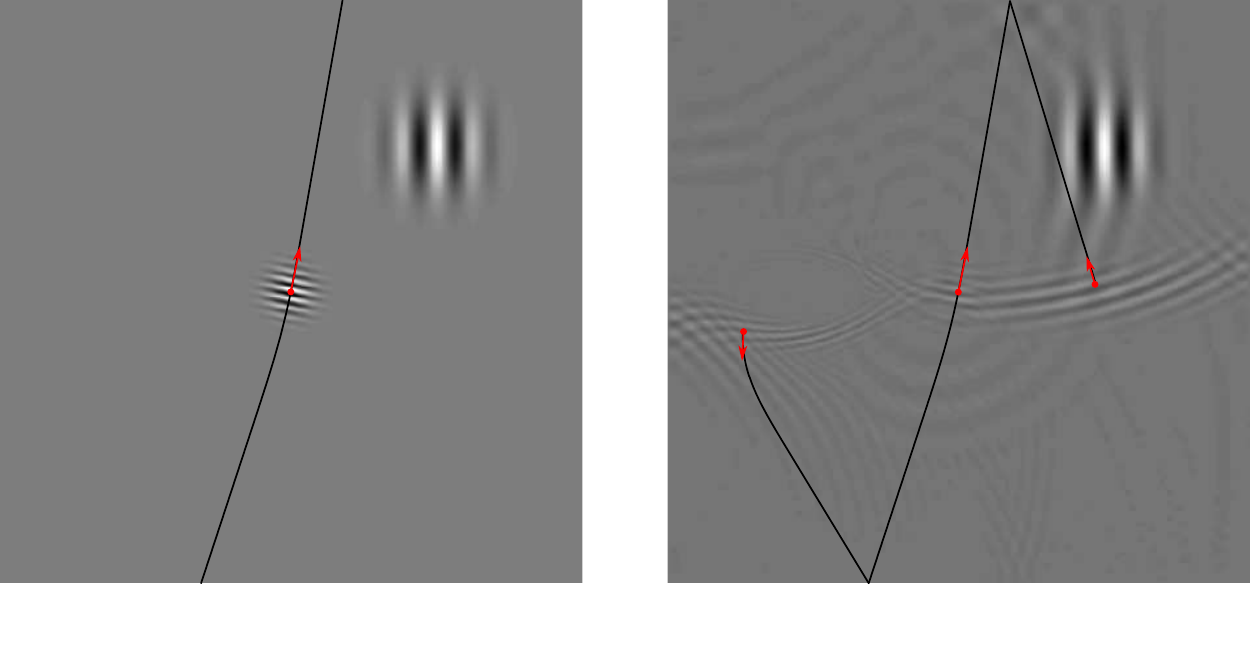\end{center}
\caption{Artifacts in a reconstructed image with $M f$ under sampled
  in time variable and a variable wave speed.  We trace the geodesics
  to find the image of $(x,\xi)$ under the map
  $C_{\pm}\circ S_k \circ C_{\pm}$ as explained above.}
\label{fig:sing_var_t}
\end{figure}
\begin{figure}[!ht]
\centering
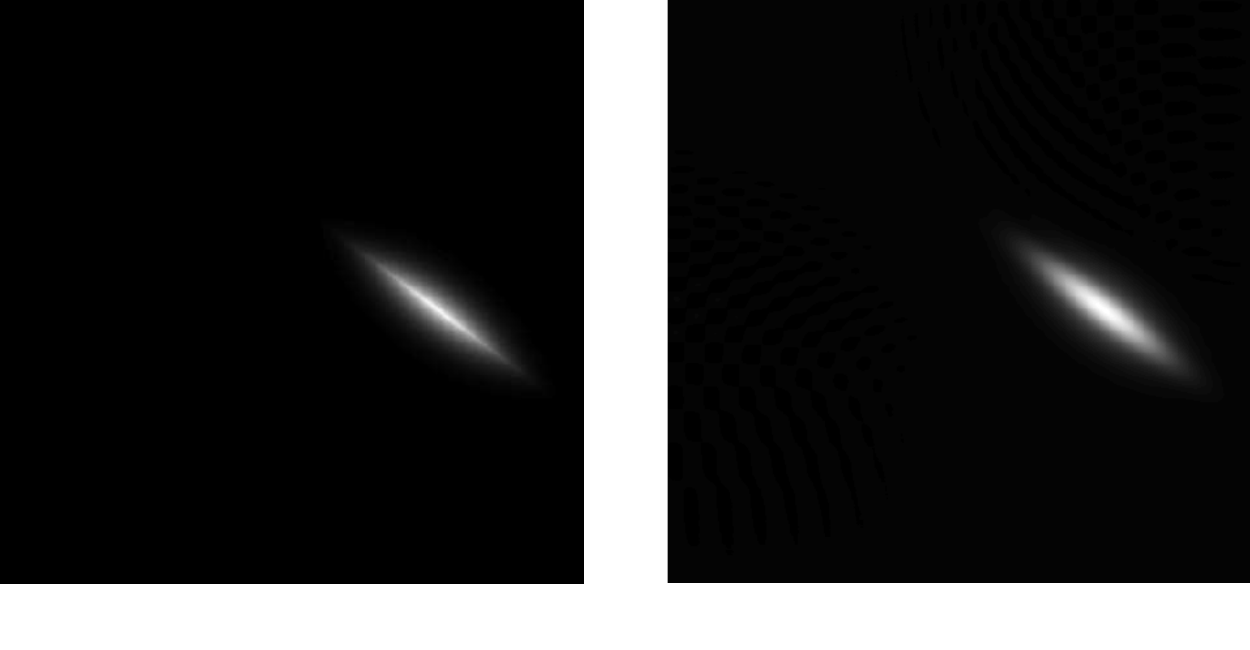
\caption{Original and reconstructed image of a smooth approximation of
  an line segment. Here we have under sampled in $t$.  The under sampling has
  resulted in blurring of this ``line segment''.  This is due to the fact that
  under sampling in $t$ shifts high frequency data in
  $\mathcal{F}(Mf)$.}
\label{fig:wormReconT}
\end{figure}
\begin{figure}[!ht]
\centering
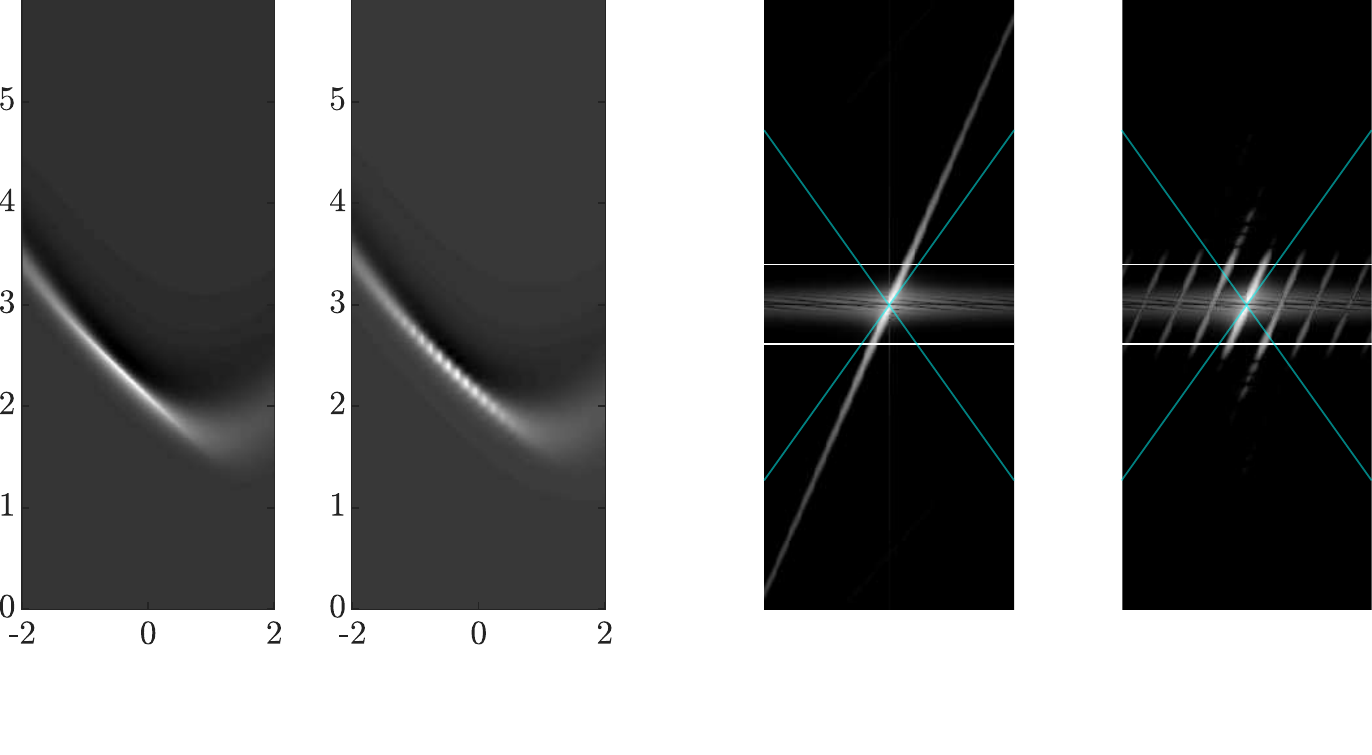
\caption{Collected data and Fourier transform along with under sampled
  data in $t$ for example given in Figure \ref{fig:wormReconT}.  Data
  was collected on all edges of the square at a rate guaranteeing no
  aliasing.  Shown is the data from the bottom edge of the square.  We
  can see that under sampling in $t$ has resulted in the Fourier
  Transform of $M f$ being folded into the band limit region.
  Under sampling in $t$ shifts large frequencies from
  $\mathcal{F}(M f)$, thus producing the blurred image we see in the
  right of Figure \ref{fig:wormReconT}.}
\label{fig:wormTData}
\end{figure}
\begin{figure}[p]
\centering
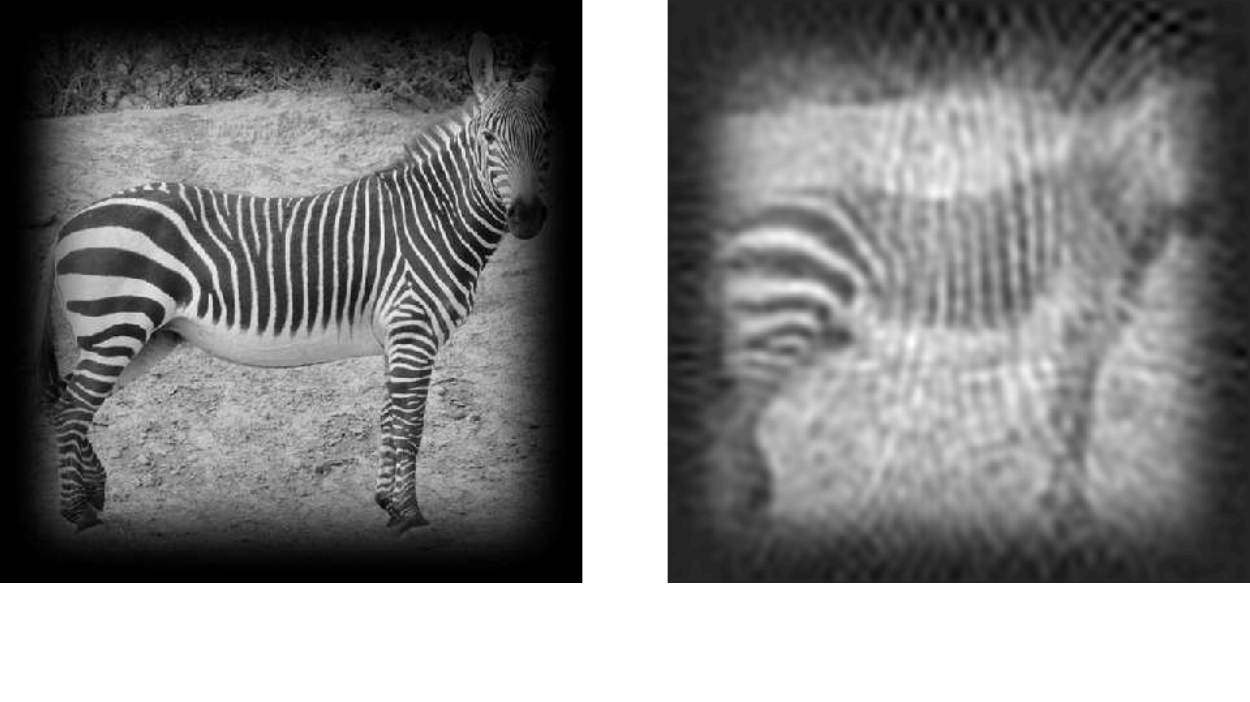
\caption{Image of a zebra along with reconstruction from under sampled
  (in $t$) data.  The wave speed here is constant.  High frequencies
  are lost due to this under sampling and the result is a heavily
  blurred image with aliasing artifacts.}
\label{fig:zebraT}
\end{figure}
\begin{figure}[p]
\centering
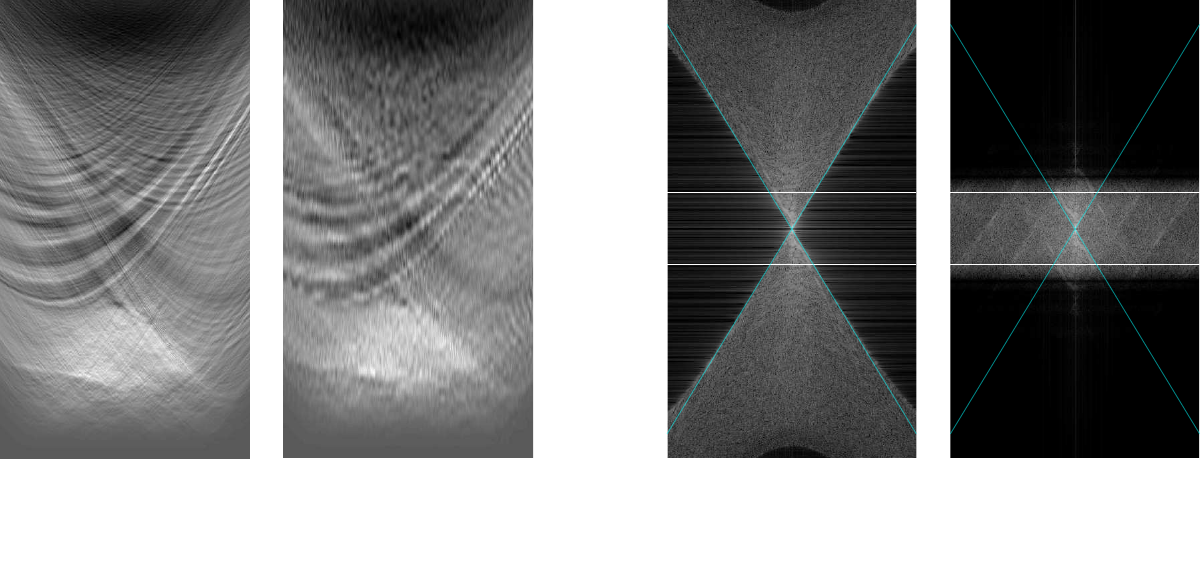
\caption{Collected data and under sampled data in $t$ along with the
  associated Fourier transform data for the zebra image above.  Note
  that the high frequencies in $\mathcal{F}(Mf)$ have be shifted so
  that they are approximately in the band $-\pi/s_t < \tau < \pi/s_t$,
  which is what results in the blurring in the reconstruction.}
\label{fig:zebraTData}
\end{figure}

We include a more complicated image reconstruction in Figure
\ref{fig:zebraT} along with the collected data in Figure
\ref{fig:zebraTData}. We also show how a smooth approximation of an
line segment is affected by these artifacts in the image given in
Figure \ref{fig:wormReconT}.  For this image and reconstruction, we
have included the collected data and Fourier transform images in
Figure \ref{fig:wormTData}.

\subsection{Under sampling in \(y\)}
\begin{figure}[!ht]
\centering
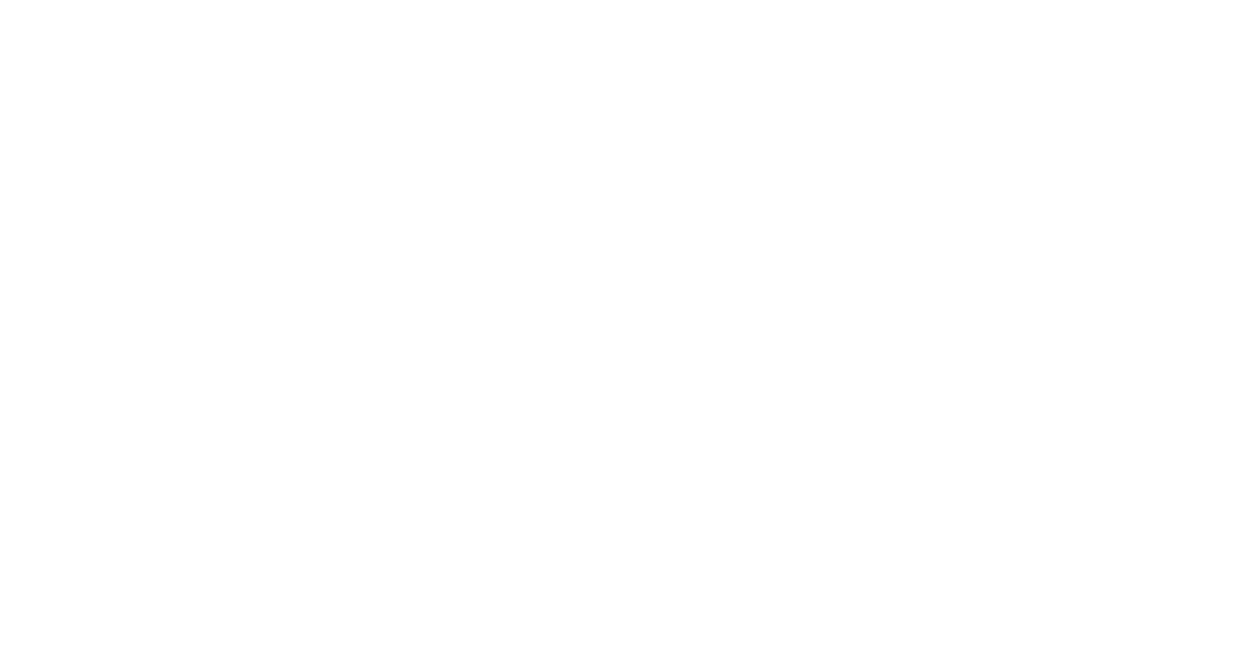
\caption{Artifacts in a reconstructed image with $M f$ under sampled
  in space variables. Here we take $c\equiv 1$. Specifically, $M f$
  here was under sampled on the left and right edges of the square.
  Note that there is no artifact in the reconstructed image coming
  from the pattern in the upper right corner of the square, because
  singularities from this pattern hit the boundary of the square
  perpendicularly. Note also that the original singularity still
  remains with half its amplitude because we did not under sample along
  the bottom edge of the square.}
\label{fig:sing_const_y}
\end{figure}

\begin{figure}[!ht]
\centering
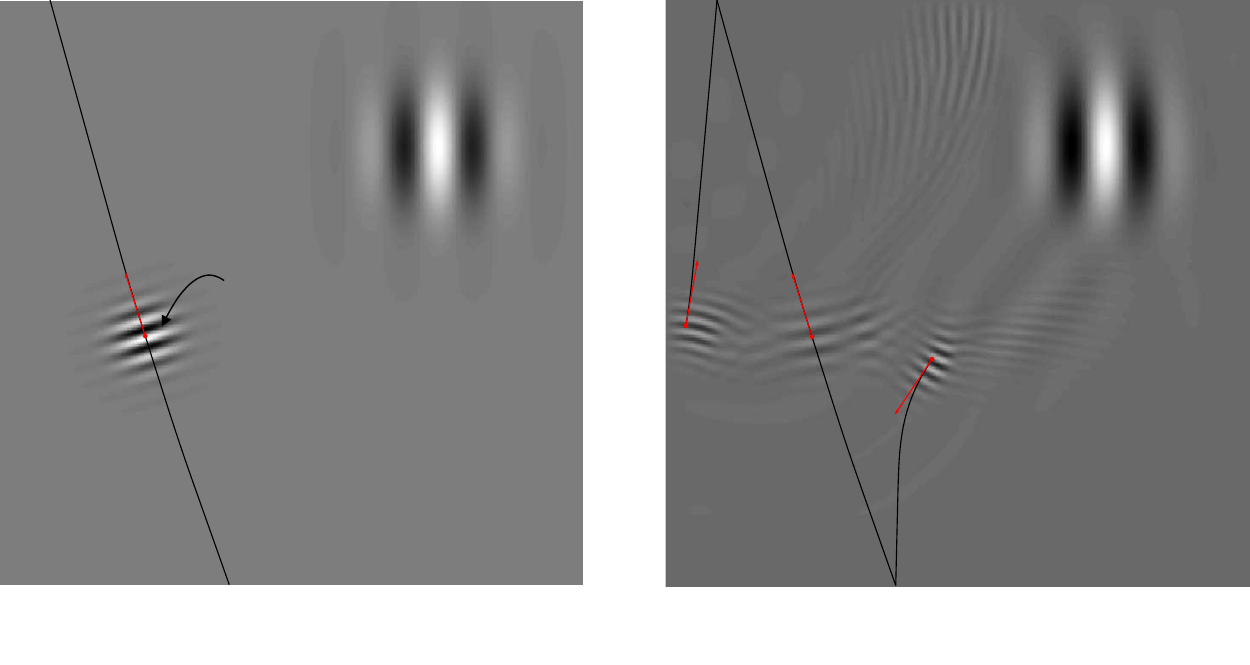
\caption{Artifacts in a reconstructed image with $M f$ under sampled
  in space variables and a variable wave speed.  Specifically, $M f$
  here was under sampled on the top and bottom edges of the square.
  The artifacts in the reconstruction have the same frequency as the
  original, but with a space shift due to under sampling.}
\label{fig:sing_var_y}
\end{figure}

Now suppose that we have under sampled the $y$ variable, i.e. we have
chosen $s_{y^j} > \frac{\pi}{B}$ for some $j = 1,\ldots, n$.  Then
again, we will have aliasing and the error in the reconstruction will
involve the frequency shift operator, but now $S_k$ will act on
$\eta_j$ as
\begin{gather*}
S_k : \eta_j \mapsto \eta_j + \frac{2\pi k}{s_{y^j}}.
\end{gather*}
This operator is valid as long as $\eta_j + \frac{2\pi k}{s_{y^j}} \in
[-\pi/s_{y^j},\pi/s_{y^j}]$.  The canonical relation of the $h$-FIO that operates on $M f$ as a
reconstruction of $f$ will then be given by (again, we only consider
$C_+$ here)
\begin{gather*}
  C_+^{-1}\circ S_k \circ C_+ (x,\xi) =
  C_+^{-1}(s_+(x,\xi),\gamma_{x,\xi}(s_+(x,\xi)),-|\xi|_g,\dot{\gamma}_{x,\xi}'(s_+(x,\xi))+\frac{2\pi
    k}{s_{y^j}}\mathbf{e}_j),
\end{gather*}
where $\mathbf{e}_j$ is the unit vector in the $y^j$ direction. Note
that, in particular, this implies that the artifacts will have the
same frequency as that of the original image, but perhaps with a space
shift.  Also, because this operator is valid as long as
$\eta_j + 2\pi k/s_{y^j} \in [-\pi/s_{y^j} , \pi/s_{y^j}]$, if the
geodesic emanating from $(x,\xi)$ hits the boundary $\partial \Omega$
perpendicularly, then the point $(x,\xi)$ will be unaffected by this
shift in the reconstruction, i.e. there will be no artifacts that come
from $(x,\xi)$.  This is true because if the geodesic emanating from
$(x,\xi)$ hits $\partial\Omega$ perpendicularly, then $\eta_j = 0$ and
$2\pi k/s_{y^j} \not\in [-\pi/s_{y^j},\pi/s_{y^j}]$ for any
$k \neq 0$.  Finding these artifacts in practice follows in much the
same way as finding where artifacts occur for under sampling in the
time variable.  We illustrate this for the constant speed, Euclidean
case in Figure \ref{fig:sing_const_y} and see Figure
\ref{fig:sing_var_y} for the variable speed case.


\begin{figure}[!ht]
\centering
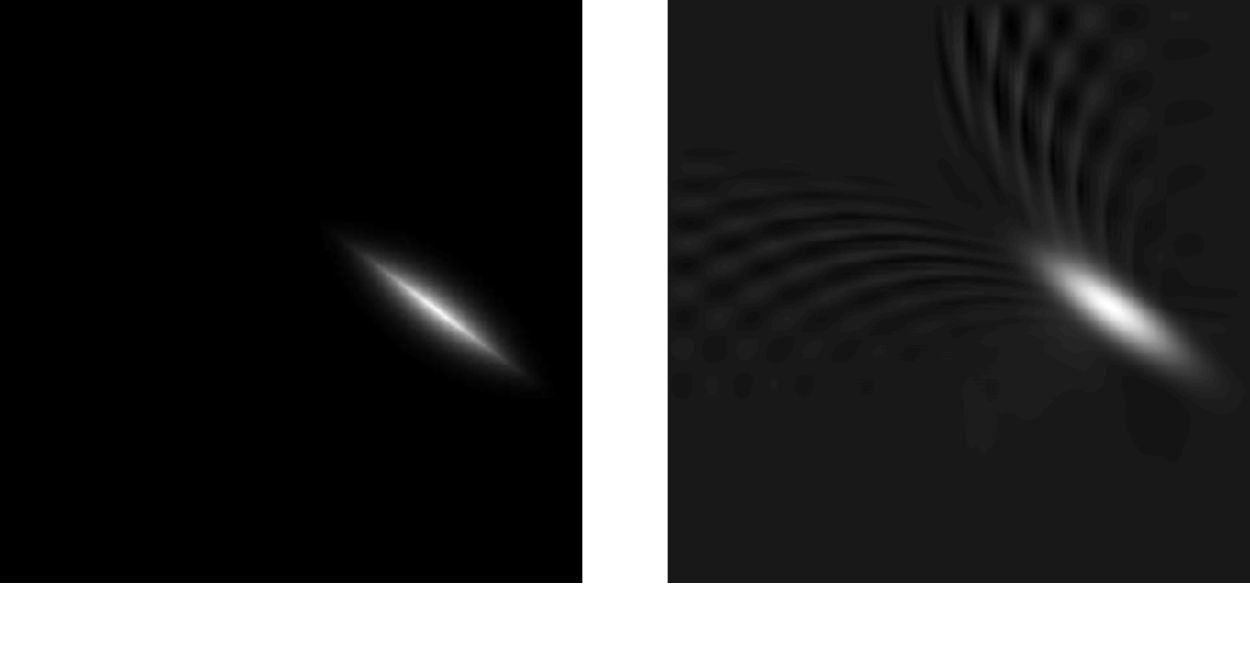
\caption{Original and reconstructed image of a smooth approximation of
  an line segment.  Here we have under sampled in $y$.  This has resulted in
  some blurring, but also in high frequency artifacts.}
\label{fig:wormReconY}
\end{figure}

\begin{figure}[!ht]
\centering
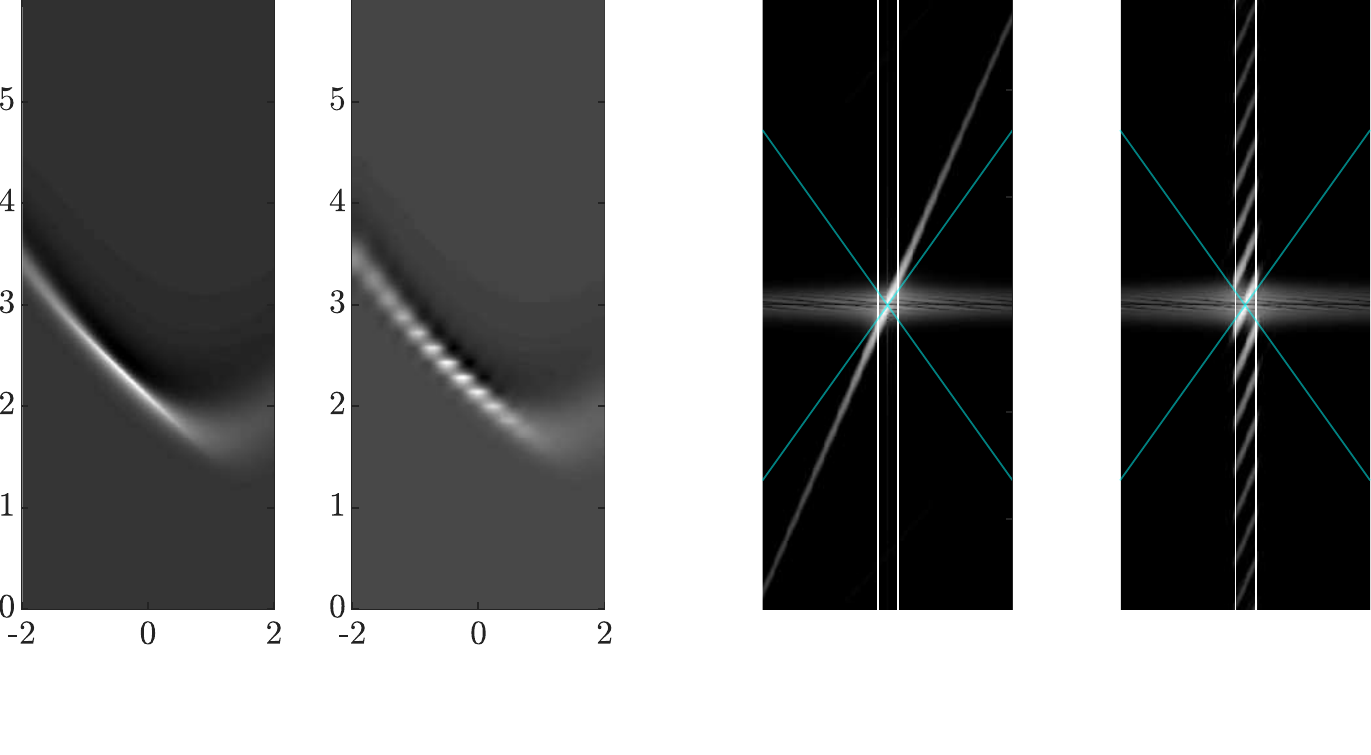
\caption{Collected data and Fourier transform along with under sampled
  data in $y$ for example given in Figure \ref{fig:wormReconY}.  In
  contrast to when we under sample in $t$, we see that high frequencies
  in $\mathcal{F}(M f)$ are not necessarily eliminated when we
  under sample in $y$, but there is a phase shift.  This results in
  more high frequency artifacts in the image on the right in Figure
  \ref{fig:wormReconY}.}
\label{fig:wormYData}
\end{figure}

\begin{figure}[p]
\centering
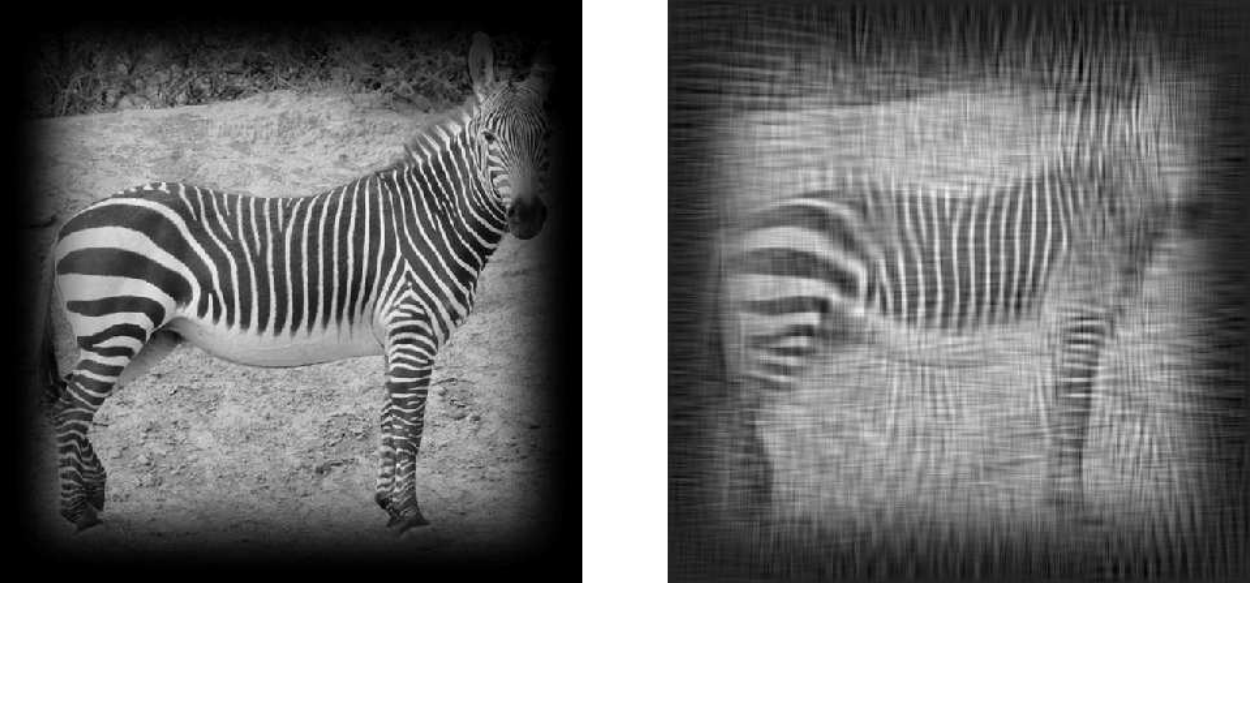
\caption{Image of a zebra along with reconstruction from under sampled
  (in $y$) data.  The wave speed here is constant.  Note that the
  singularities that hit the boundary of the square nearly
  perpendicularly are preserved, but there are also a lot of high
  frequency artifacts in the reconstructed image.}
\label{fig:zebraY}
\end{figure}

\begin{figure}[p]
\centering
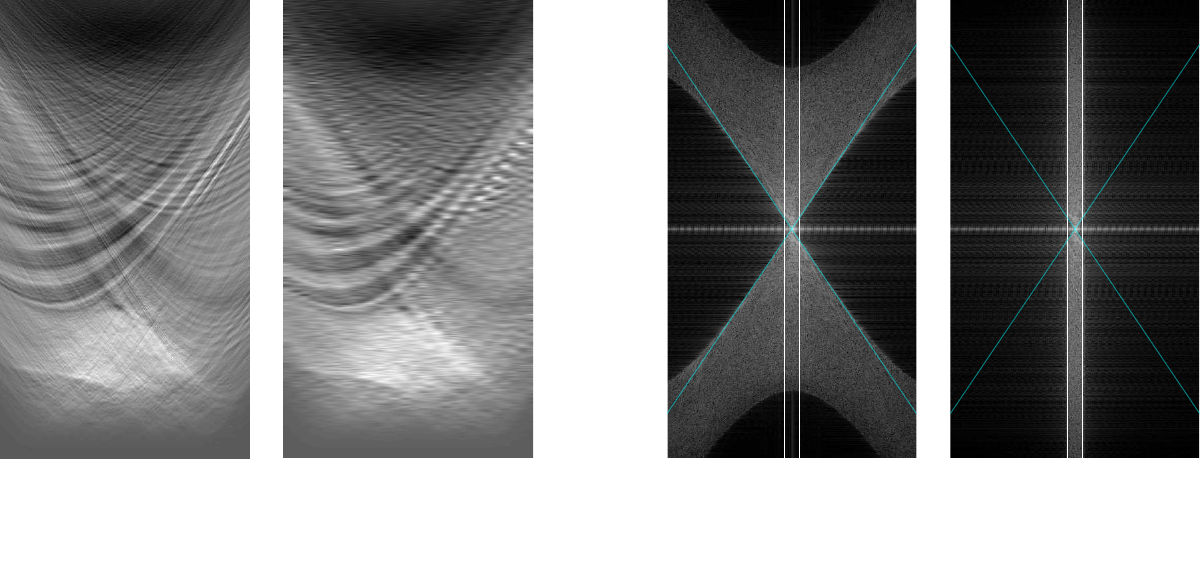
\caption{Collected data and under sampled data in $y$ along with the
  associated Fourier transform data for the zebra image above.
  Under sampling has resulted in the shifting of frequencies in
  $\mathcal{F}(Mf)$ so that $-\pi/s_{y^j} < \eta < \pi/s_{y^j}$.  This moves
  high frequencies but does not destroy them, which is what causes the
  high frequency artifacts in the reconstructed image above.}
\label{fig:zebraYData}
\end{figure}

We again include a more complicated image reconstruction in Figure
\ref{fig:zebraY} along with the collected data in Figure
\ref{fig:zebraYData}.  We also show how a smooth approximation of an
line segment is affected by these artifacts in the image given in
Figure \ref{fig:wormReconY}.  For this image and reconstruction, we
have included the collected data and Fourier transform images in
Figure \ref{fig:wormYData}.

\section{Averaged data}
\begin{figure}[!ht]
\centering
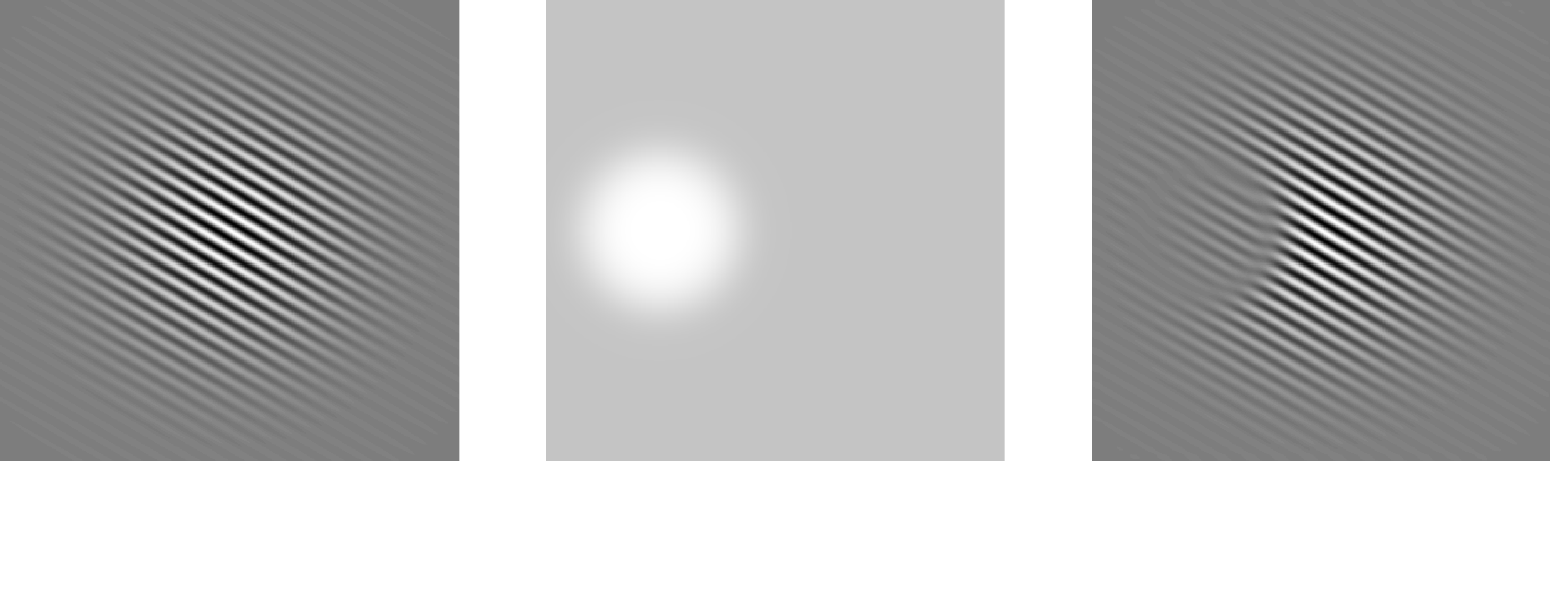
\caption{Reconstructed image from data that has been averaged in time
  variable.  We can see that the reconstructed image is most blurred
  at the points where the speed $c(x)$ is fast, and there is less
  blurring where $c(x) = 1$.}
\label{fig:avg_var_t}
\end{figure}

Suppose that the collected data $M f(t,y)$ has been averaged in the
$t$ or $y$ variables for some reason (in practice this can be done to
try to avoid aliasing, or in an attempt to reduce the noise
in data).  This can be modeled in a few ways, including taking a
convolution $\phi_h * M f$ with a smooth function
$\phi_h = h^n \phi(\cdot / h)$ that decreases away from the origin to
0.  To model localized averaging however, we will consider data of the
form $Q_h M f (t,y)$, where $Q_h$ is an $h$-$\Psi$DO with a principal
symbol of the form
$q_0(t,y,\tau,\eta) = \psi(a |\tau|^2 + b |\eta|^2)$ where
$\psi \in C_0^\infty( \R)$ is decreasing.  The effect of $Q_h$ is to
limit $\WF_h (M f)$, which will in principle remove the high frequency
singularities of $M f$ which will have a smoothing effect.  From
\cite{StefanovP2018}, we know that because $M$ is a FIO associated
with the canonical map $C = C_+ \cup C_-$, that the composition
$Q_hM f$ can be written
\begin{gather*}
Q_hM f = M P_h f + O(h^\infty) f,
\end{gather*}
where $P_h$ is a $h$-$\Psi$DO with principal symbol
$p_0 = q_0 \circ C$ where $q_0$ is the principal symbol of $Q_h$.  So,
for $Q_h, q_0$, we may calculate
\begin{align*}
  p_0(x,\xi) &= \frac{1}{2}\left( q_0\circ C_+(x,\xi) + q_0 \circ C_-(x,\xi)\right) \\
             &= \frac{1}{2} \left( \psi(a |\xi|_g^2 + b|\dot{\gamma}_{x,\xi}'(s_+(x,\xi))|_{g_{0,\partial\Omega}}^2) + \psi(a |\xi|_g^2 + b|\dot{\gamma}_{x,\xi}'(s_-(x,\xi))|_{g_{0,\partial\Omega}}^2)\right).
\end{align*}
Suppose we only average the time data in $M f (t,y)$.  This
corresponds to taking $b = 0$ above to give
$p_0(x,\xi) = \psi(a |\xi|_g^2)$.  This symbol takes its minimum
values where $|\xi|_g^2 = c^2(x)g_0^{ij}\xi_i\xi_j$ is maximized.
Assuming for a moment that $g$ is Euclidean, this means that we expect
more blurring at points $(x,\xi)$ where the wave speed is ``fast''.
Additionally, we expect singularities $(x,\xi)$ with large frequencies
$|\xi|$ to be blurred more than smaller frequencies where the wave
speed is the same.  These effects can both be seen in Figure
\ref{fig:avg_var_t}.

Suppose now that we only average data in the spatial variable $y$.
This corresponds to taking $a=0$ above and we get the principle symbol
of $p_0$ to be
$$
p_0(x,\xi) = \frac{1}{2}\left( \psi(b
  |\dot{\gamma}_{x,\xi}'(s_+(x,\xi))|^2) + \psi(b
  |\dot{\gamma}_{x,\xi}'(s_-(x,\xi))|^2) \right).
$$
Here the norm is the induced norm on the boundary, which we have noted
in this paper as $g_{0,\partial\Omega}$.  This symbol takes its
smallest values when $|\dot{\gamma}_{x,\xi}'(s_\pm(x,\xi))|^2$
is large, i.e. when the geodesic issued from $(x,\xi)$ intersects the
boundary $\partial \Omega$ at a large angle.  In addition, we expect
singularities that hit the boundary $\partial \Omega$ perpendicularly
to be affected far less by averaging of data in the $y$ variable.  In
addition, because
$|\dot{\gamma}'_{x,\xi}(s_{\pm}(x,\xi))|^2 = |\xi|_g^2
\cos^2(\theta_{\pm})$ where $\theta_{\pm}$ is the angle between
$\dot{\gamma}'_{x,\xi}(s_{\pm}(x,\xi))$ and
$\dot{\gamma}_{x,\xi}(s_{\pm}(x,\xi))$ we expect to see more blurring
at points with faster speeds or higher frequency.  For constant speeds
$c$, the effect of averaging data in $t$ is uniform in $\Omega$, but
the effect is local for averaging in $y$, due to the blurring
depending on the angle of intersection made by geodesics.  In
addition, with a variable speed singularities in ``slow spots'' of $c$
will have higher resolution when blurring $Mf(t,y)$ in the $y$-data,
but their resolution will still depend on how geodesics hit the
boundary.  The result is a roughly uniform blurring in fast spots of
$c$, and local blurring elsewhere in the image depending on the
geometry determined by $c^{-2}g_0$. This can be seen in Figure
\ref{fig:avg_var_y} below.

\begin{figure}[!ht]
\centering
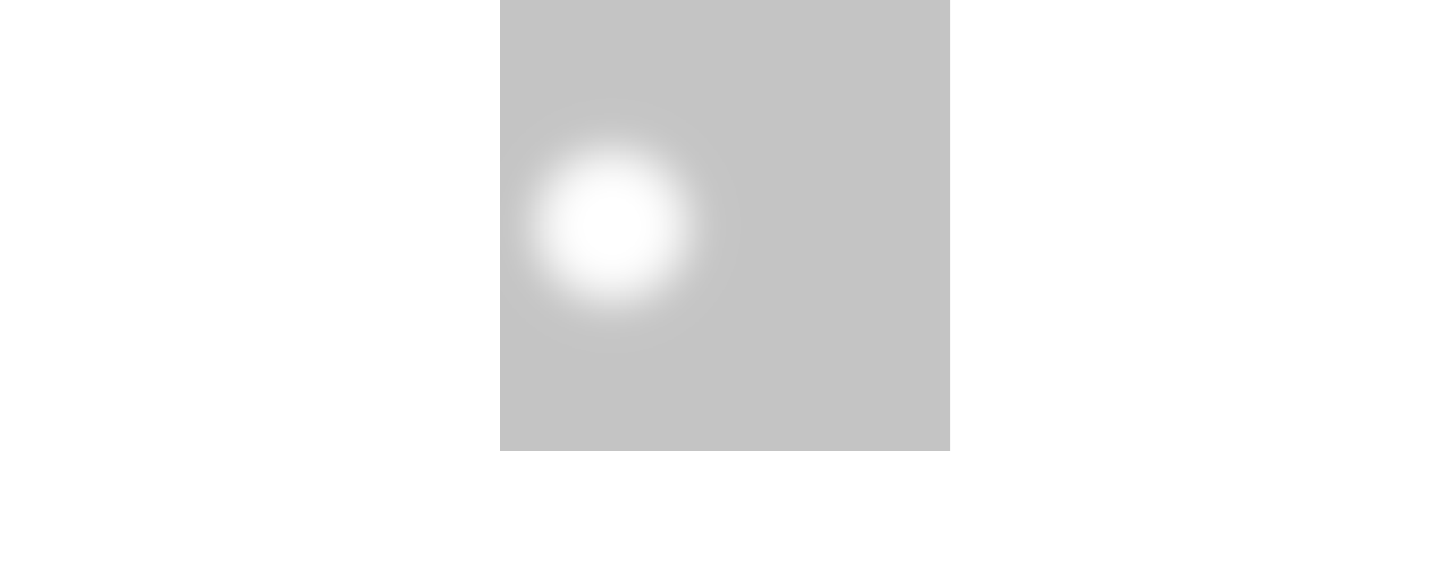
\caption{Reconstructed image from data that has been averaged in space
  variable.  We can see from the drawn in geodesics, that
  singularities that hit the boundary at a larger angle to the normal
  vector to the boundary are blurred more in the reconstructed image
  after averaging the collected data.  Meanwhile, those singularities
  that hit the boundary nearly perpendicularly are largely unaffected
  by the averaging of the data, at least on one side.}
\label{fig:avg_var_y}
\end{figure}

\section{Anti-aliasing}
\begin{figure}[!ht]
\centering
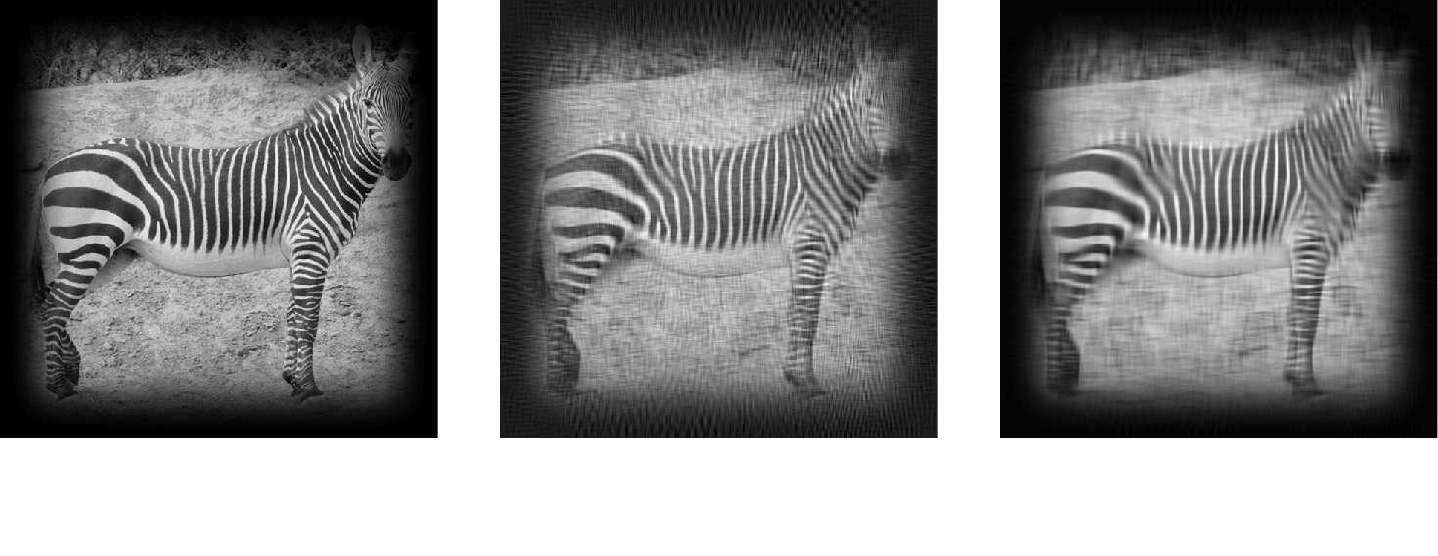
\caption{Example showing anti-aliasing scheme in which we first
  average the data $Mf(t,y)$ in the $y$ variable and then sample this
  blurred version given by $Q_h Mf (t,y)$ in the above notation.  We
  can see that some of the aliasing artifacts have been
  removed at the cost of some loss of resolution.}
\label{fig:antiAlias}
\end{figure}

We can use the above discussion to propose an anti-aliasing scheme.
Averaging the measured data $Mf(t,y)$ in the space variable can be
accomplished in practice in many ways, whether by using small
averaging detectors, or by vibrating the boundary $\partial \Omega$
where we are taking pointwise measurements.  We know then that this
can be modeled by applying the $h$-$\Psi$DO $Q_h$ to $Mf$ where is as
in the previous section.  This then allows us to say that
$Q_h Mf(t,y) = MP_h f(t,y) + O(h^\infty)f$. In other words, by
averaging the data in $y$, we measure $P_h f (x,\xi)$, where $P_h$ is
an $h$-$\Psi$DO with principle symbol
$p_0(x,\xi) = q_0 \circ C(x,\xi)$ and $C$ is the canonical relation of
$M$, plus some error term with low order frequencies.  We then expect
that if we average $Mf(t,y)$ in the $y$ variable before sampling, this
should remove some of the shifting aliasing artifacts that
appear when $Mf(t,y)$ has been under sampled in $y$, perhaps at the
cost of some loss of resolution.  See Figure \ref{fig:antiAlias} for
an example of this anti-aliasing scheme in action.

As a final note, we point out that under sampling in the time variable
$t$ can cause data in $\mathcal{F}(Mf)$ to shift outside of the
characteristic cone, and from this, one should be able to recover some
high frequency singularities from data $Mf(t,y)$ under sampled in $t$
by shifting these singularities back out into the characteristic cone
where they necessarily originated (see Figure \ref{fig:charCone}).
However, we can only recover a small fraction of the high frequency
singularities in this way uniquely in special cases, and in general we
cannot recover the singularities without adding high frequency
artifacts to the reconstructed image.

\printbibliography
\end{document}

%% file: resolutionFastV2.pdf_tex
\begingroup%
  \makeatletter%
  \providecommand\color[2][]{%
    \errmessage{(Inkscape) Color is used for the text in Inkscape, but the package 'color.sty' is not loaded}%
    \renewcommand\color[2][]{}%
  }%
  \providecommand\transparent[1]{%
    \errmessage{(Inkscape) Transparency is used (non-zero) for the text in Inkscape, but the package 'transparent.sty' is not loaded}%
    \renewcommand\transparent[1]{}%
  }%
  \providecommand\rotatebox[2]{#2}%
  \newcommand*\fsize{\dimexpr\f@size pt\relax}%
  \newcommand*\lineheight[1]{\fontsize{\fsize}{#1\fsize}\selectfont}%
  \ifx\svgwidth\undefined%
    \setlength{\unitlength}{446.39998627bp}%
    \ifx\svgscale\undefined%
      \relax%
    \else%
      \setlength{\unitlength}{\unitlength * \real{\svgscale}}%
    \fi%
  \else%
    \setlength{\unitlength}{\svgwidth}%
  \fi%
  \global\let\svgwidth\undefined%
  \global\let\svgscale\undefined%
  \makeatother%
  \begin{picture}(1,0.35710809)%
    \lineheight{1}%
    \setlength\tabcolsep{0pt}%
    \put(0,0){\includegraphics[width=\unitlength,page=1]{resolutionFastV2.pdf}}%
    \put(0.06624051,0.00037638){\color[rgb]{0,0,0}\makebox(0,0)[lt]{\lineheight{1.25}\smash{\begin{tabular}[t]{l}Original Image\end{tabular}}}}%
    \put(0.33710478,0.00037638){\color[rgb]{0,0,0}\makebox(0,0)[lt]{\lineheight{1.25}\smash{\begin{tabular}[t]{l}Time Reversal Reconstruction\end{tabular}}}}%
    \put(0.7043953,0.00037638){\color[rgb]{0,0,0}\makebox(0,0)[lt]{\lineheight{1.25}\smash{\begin{tabular}[t]{l}Wave Speed with Fast Spot\end{tabular}}}}%
    \put(-0.01241214,0.03728495){\color[rgb]{0,0,0}\makebox(0,0)[lt]{\lineheight{1.25}\smash{\begin{tabular}[t]{l}-2\end{tabular}}}}%
    \put(0.28999266,0.03728495){\color[rgb]{0,0,0}\makebox(0,0)[lt]{\lineheight{1.25}\smash{\begin{tabular}[t]{l}2\end{tabular}}}}%
    \put(0.34014292,0.03728495){\color[rgb]{0,0,0}\makebox(0,0)[lt]{\lineheight{1.25}\smash{\begin{tabular}[t]{l}-2\end{tabular}}}}%
    \put(0.64254773,0.03728495){\color[rgb]{0,0,0}\makebox(0,0)[lt]{\lineheight{1.25}\smash{\begin{tabular}[t]{l}2\end{tabular}}}}%
    \put(0.69269796,0.03728495){\color[rgb]{0,0,0}\makebox(0,0)[lt]{\lineheight{1.25}\smash{\begin{tabular}[t]{l}-2\end{tabular}}}}%
    \put(0.99510278,0.03728495){\color[rgb]{0,0,0}\makebox(0,0)[lt]{\lineheight{1.25}\smash{\begin{tabular}[t]{l}2\end{tabular}}}}%
    \put(-0.02928884,0.05313341){\color[rgb]{0,0,0}\makebox(0,0)[lt]{\lineheight{1.25}\smash{\begin{tabular}[t]{l}-2\end{tabular}}}}%
    \put(-0.02031077,0.3453828){\color[rgb]{0,0,0}\makebox(0,0)[lt]{\lineheight{1.25}\smash{\begin{tabular}[t]{l}2\end{tabular}}}}%
    \put(0.32374349,0.05313341){\color[rgb]{0,0,0}\makebox(0,0)[lt]{\lineheight{1.25}\smash{\begin{tabular}[t]{l}-2\end{tabular}}}}%
    \put(0.33272156,0.3453828){\color[rgb]{0,0,0}\makebox(0,0)[lt]{\lineheight{1.25}\smash{\begin{tabular}[t]{l}2\end{tabular}}}}%
    \put(0.6712833,0.05313341){\color[rgb]{0,0,0}\makebox(0,0)[lt]{\lineheight{1.25}\smash{\begin{tabular}[t]{l}-2\end{tabular}}}}%
    \put(0.68026137,0.3453828){\color[rgb]{0,0,0}\makebox(0,0)[lt]{\lineheight{1.25}\smash{\begin{tabular}[t]{l}2\end{tabular}}}}%
  \end{picture}%
\endgroup%

%% file: resolutionSlowV2.pdf_tex
\begingroup%
  \makeatletter%
  \providecommand\color[2][]{%
    \errmessage{(Inkscape) Color is used for the text in Inkscape, but the package 'color.sty' is not loaded}%
    \renewcommand\color[2][]{}%
  }%
  \providecommand\transparent[1]{%
    \errmessage{(Inkscape) Transparency is used (non-zero) for the text in Inkscape, but the package 'transparent.sty' is not loaded}%
    \renewcommand\transparent[1]{}%
  }%
  \providecommand\rotatebox[2]{#2}%
  \newcommand*\fsize{\dimexpr\f@size pt\relax}%
  \newcommand*\lineheight[1]{\fontsize{\fsize}{#1\fsize}\selectfont}%
  \ifx\svgwidth\undefined%
    \setlength{\unitlength}{446.39998627bp}%
    \ifx\svgscale\undefined%
      \relax%
    \else%
      \setlength{\unitlength}{\unitlength * \real{\svgscale}}%
    \fi%
  \else%
    \setlength{\unitlength}{\svgwidth}%
  \fi%
  \global\let\svgwidth\undefined%
  \global\let\svgscale\undefined%
  \makeatother%
  \begin{picture}(1,0.35710809)%
    \lineheight{1}%
    \setlength\tabcolsep{0pt}%
    \put(0,0){\includegraphics[width=\unitlength,page=1]{resolutionSlowV2.pdf}}%
    \put(0.0662405,0.00000003){\color[rgb]{0,0,0}\makebox(0,0)[lt]{\lineheight{1.25}\smash{\begin{tabular}[t]{l}Original Image\end{tabular}}}}%
    \put(0.33710483,0.00000003){\color[rgb]{0,0,0}\makebox(0,0)[lt]{\lineheight{1.25}\smash{\begin{tabular}[t]{l}Time Reversal Reconstruction\end{tabular}}}}%
    \put(0.69917125,0.00000003){\color[rgb]{0,0,0}\makebox(0,0)[lt]{\lineheight{1.25}\smash{\begin{tabular}[t]{l}Wave Speed with Slow Spot\end{tabular}}}}%
    \put(-0.01241261,0.03728499){\color[rgb]{0,0,0}\makebox(0,0)[lt]{\lineheight{1.25}\smash{\begin{tabular}[t]{l}-2\end{tabular}}}}%
    \put(0.28999217,0.03728499){\color[rgb]{0,0,0}\makebox(0,0)[lt]{\lineheight{1.25}\smash{\begin{tabular}[t]{l}2\end{tabular}}}}%
    \put(0.34014243,0.03728499){\color[rgb]{0,0,0}\makebox(0,0)[lt]{\lineheight{1.25}\smash{\begin{tabular}[t]{l}-2\end{tabular}}}}%
    \put(0.64254723,0.03728499){\color[rgb]{0,0,0}\makebox(0,0)[lt]{\lineheight{1.25}\smash{\begin{tabular}[t]{l}2\end{tabular}}}}%
    \put(0.69269746,0.03728499){\color[rgb]{0,0,0}\makebox(0,0)[lt]{\lineheight{1.25}\smash{\begin{tabular}[t]{l}-2\end{tabular}}}}%
    \put(0.99510225,0.03728499){\color[rgb]{0,0,0}\makebox(0,0)[lt]{\lineheight{1.25}\smash{\begin{tabular}[t]{l}2\end{tabular}}}}%
    \put(-0.02928931,0.05313345){\color[rgb]{0,0,0}\makebox(0,0)[lt]{\lineheight{1.25}\smash{\begin{tabular}[t]{l}-2\end{tabular}}}}%
    \put(-0.02031124,0.34538283){\color[rgb]{0,0,0}\makebox(0,0)[lt]{\lineheight{1.25}\smash{\begin{tabular}[t]{l}2\end{tabular}}}}%
    \put(0.323743,0.05313345){\color[rgb]{0,0,0}\makebox(0,0)[lt]{\lineheight{1.25}\smash{\begin{tabular}[t]{l}-2\end{tabular}}}}%
    \put(0.33272107,0.34538283){\color[rgb]{0,0,0}\makebox(0,0)[lt]{\lineheight{1.25}\smash{\begin{tabular}[t]{l}2\end{tabular}}}}%
    \put(0.67128282,0.05313345){\color[rgb]{0,0,0}\makebox(0,0)[lt]{\lineheight{1.25}\smash{\begin{tabular}[t]{l}-2\end{tabular}}}}%
    \put(0.68026089,0.34538283){\color[rgb]{0,0,0}\makebox(0,0)[lt]{\lineheight{1.25}\smash{\begin{tabular}[t]{l}2\end{tabular}}}}%
  \end{picture}%
\endgroup%

%% file: resolutionFastY.pdf_tex
\begingroup%
  \makeatletter%
  \providecommand\color[2][]{%
    \errmessage{(Inkscape) Color is used for the text in Inkscape, but the package 'color.sty' is not loaded}%
    \renewcommand\color[2][]{}%
  }%
  \providecommand\transparent[1]{%
    \errmessage{(Inkscape) Transparency is used (non-zero) for the text in Inkscape, but the package 'transparent.sty' is not loaded}%
    \renewcommand\transparent[1]{}%
  }%
  \providecommand\rotatebox[2]{#2}%
  \newcommand*\fsize{\dimexpr\f@size pt\relax}%
  \newcommand*\lineheight[1]{\fontsize{\fsize}{#1\fsize}\selectfont}%
  \ifx\svgwidth\undefined%
    \setlength{\unitlength}{446.39998627bp}%
    \ifx\svgscale\undefined%
      \relax%
    \else%
      \setlength{\unitlength}{\unitlength * \real{\svgscale}}%
    \fi%
  \else%
    \setlength{\unitlength}{\svgwidth}%
  \fi%
  \global\let\svgwidth\undefined%
  \global\let\svgscale\undefined%
  \makeatother%
  \begin{picture}(1,0.35710809)%
    \lineheight{1}%
    \setlength\tabcolsep{0pt}%
    \put(0.06682375,0.0003764){\color[rgb]{0,0,0}\makebox(0,0)[lt]{\lineheight{1.25}\smash{\begin{tabular}[t]{l}Original Image\end{tabular}}}}%
    \put(0.33804187,0.0003764){\color[rgb]{0,0,0}\makebox(0,0)[lt]{\lineheight{1.25}\smash{\begin{tabular}[t]{l}Time Reversal Reconstruction\end{tabular}}}}%
    \put(0.70450035,0.0003764){\color[rgb]{0,0,0}\makebox(0,0)[lt]{\lineheight{1.25}\smash{\begin{tabular}[t]{l}Wave Speed with Fast Spot\end{tabular}}}}%
    \put(0,0){\includegraphics[width=\unitlength,page=1]{resolutionFastY.pdf}}%
    \put(-0.01241261,0.03728495){\color[rgb]{0,0,0}\makebox(0,0)[lt]{\lineheight{1.25}\smash{\begin{tabular}[t]{l}-2\end{tabular}}}}%
    \put(0.28999217,0.03728495){\color[rgb]{0,0,0}\makebox(0,0)[lt]{\lineheight{1.25}\smash{\begin{tabular}[t]{l}2\end{tabular}}}}%
    \put(0.34014246,0.03728495){\color[rgb]{0,0,0}\makebox(0,0)[lt]{\lineheight{1.25}\smash{\begin{tabular}[t]{l}-2\end{tabular}}}}%
    \put(0.64254726,0.03728495){\color[rgb]{0,0,0}\makebox(0,0)[lt]{\lineheight{1.25}\smash{\begin{tabular}[t]{l}2\end{tabular}}}}%
    \put(0.69269749,0.03728495){\color[rgb]{0,0,0}\makebox(0,0)[lt]{\lineheight{1.25}\smash{\begin{tabular}[t]{l}-2\end{tabular}}}}%
    \put(0.9951023,0.03728495){\color[rgb]{0,0,0}\makebox(0,0)[lt]{\lineheight{1.25}\smash{\begin{tabular}[t]{l}2\end{tabular}}}}%
    \put(-0.02928932,0.05313341){\color[rgb]{0,0,0}\makebox(0,0)[lt]{\lineheight{1.25}\smash{\begin{tabular}[t]{l}-2\end{tabular}}}}%
    \put(-0.02031124,0.34538279){\color[rgb]{0,0,0}\makebox(0,0)[lt]{\lineheight{1.25}\smash{\begin{tabular}[t]{l}2\end{tabular}}}}%
    \put(0.32374303,0.05313341){\color[rgb]{0,0,0}\makebox(0,0)[lt]{\lineheight{1.25}\smash{\begin{tabular}[t]{l}-2\end{tabular}}}}%
    \put(0.3327211,0.34538279){\color[rgb]{0,0,0}\makebox(0,0)[lt]{\lineheight{1.25}\smash{\begin{tabular}[t]{l}2\end{tabular}}}}%
    \put(0.67128283,0.05313341){\color[rgb]{0,0,0}\makebox(0,0)[lt]{\lineheight{1.25}\smash{\begin{tabular}[t]{l}-2\end{tabular}}}}%
    \put(0.6802609,0.34538279){\color[rgb]{0,0,0}\makebox(0,0)[lt]{\lineheight{1.25}\smash{\begin{tabular}[t]{l}2\end{tabular}}}}%
  \end{picture}%
\endgroup%

%% file: char_cone.pdf_tex
\begingroup%
  \makeatletter%
  \providecommand\color[2][]{%
    \errmessage{(Inkscape) Color is used for the text in Inkscape, but the package 'color.sty' is not loaded}%
    \renewcommand\color[2][]{}%
  }%
  \providecommand\transparent[1]{%
    \errmessage{(Inkscape) Transparency is used (non-zero) for the text in Inkscape, but the package 'transparent.sty' is not loaded}%
    \renewcommand\transparent[1]{}%
  }%
  \providecommand\rotatebox[2]{#2}%
  \newcommand*\fsize{\dimexpr\f@size pt\relax}%
  \newcommand*\lineheight[1]{\fontsize{\fsize}{#1\fsize}\selectfont}%
  \ifx\svgwidth\undefined%
    \setlength{\unitlength}{308.35621271bp}%
    \ifx\svgscale\undefined%
      \relax%
    \else%
      \setlength{\unitlength}{\unitlength * \real{\svgscale}}%
    \fi%
  \else%
    \setlength{\unitlength}{\svgwidth}%
  \fi%
  \global\let\svgwidth\undefined%
  \global\let\svgscale\undefined%
  \makeatother%
  \begin{picture}(1,0.46675335)%
    \lineheight{1}%
    \setlength\tabcolsep{0pt}%
    \put(0,0){\includegraphics[width=\unitlength,page=1]{char_cone.pdf}}%
    \put(0.39972669,0.20776471){\color[rgb]{0,0,0}\makebox(0,0)[lt]{\lineheight{1.25}\smash{\begin{tabular}[t]{l}$\eta$\end{tabular}}}}%
    \put(0.19675326,0.42740514){\color[rgb]{0,0,0}\makebox(0,0)[t]{\lineheight{1.25}\smash{\begin{tabular}[t]{c}$\tau$\end{tabular}}}}%
    \put(0.91153858,0.20776471){\color[rgb]{0,0,0}\makebox(0,0)[lt]{\lineheight{1.25}\smash{\begin{tabular}[t]{l}$\eta$\end{tabular}}}}%
    \put(0.71069376,0.42740514){\color[rgb]{0,0,0}\makebox(0,0)[t]{\lineheight{1.25}\smash{\begin{tabular}[t]{c}$\tau$\end{tabular}}}}%
    \put(0.71471026,0.38963803){\color[rgb]{0,0,0}\makebox(0,0)[lt]{\lineheight{1.25}\smash{\begin{tabular}[t]{l}$B$\end{tabular}}}}%
    \put(0,0){\includegraphics[width=\unitlength,page=2]{char_cone.pdf}}%
    \put(0.44347858,0.40048209){\color[rgb]{0,0,0}\makebox(0,0)[lt]{\lineheight{1.25}\smash{\begin{tabular}[t]{l}$S_k$\end{tabular}}}}%
    \put(0.92135457,0.04545348){\color[rgb]{0,0,0}\makebox(0,0)[lt]{\lineheight{1.25}\smash{\begin{tabular}[t]{l}-$\frac{\pi}{s_t}$\end{tabular}}}}%
    \put(0.92116001,0.37836961){\color[rgb]{0,0,0}\makebox(0,0)[lt]{\lineheight{1.25}\smash{\begin{tabular}[t]{l}$\frac{\pi}{s_t}$\end{tabular}}}}%
  \end{picture}%
\endgroup%

%% file: geodesics.pdf_tex
\begingroup%
  \makeatletter%
  \providecommand\color[2][]{%
    \errmessage{(Inkscape) Color is used for the text in Inkscape, but the package 'color.sty' is not loaded}%
    \renewcommand\color[2][]{}%
  }%
  \providecommand\transparent[1]{%
    \errmessage{(Inkscape) Transparency is used (non-zero) for the text in Inkscape, but the package 'transparent.sty' is not loaded}%
    \renewcommand\transparent[1]{}%
  }%
  \providecommand\rotatebox[2]{#2}%
  \newcommand*\fsize{\dimexpr\f@size pt\relax}%
  \newcommand*\lineheight[1]{\fontsize{\fsize}{#1\fsize}\selectfont}%
  \ifx\svgwidth\undefined%
    \setlength{\unitlength}{360bp}%
    \ifx\svgscale\undefined%
      \relax%
    \else%
      \setlength{\unitlength}{\unitlength * \real{\svgscale}}%
    \fi%
  \else%
    \setlength{\unitlength}{\svgwidth}%
  \fi%
  \global\let\svgwidth\undefined%
  \global\let\svgscale\undefined%
  \makeatother%
  \begin{picture}(1,0.52972846)%
    \lineheight{1}%
    \setlength\tabcolsep{0pt}%
    \put(0,0){\includegraphics[width=\unitlength,page=1]{geodesics.pdf}}%
    \put(0.59493804,0.47584762){\color[rgb]{0,0,0}\makebox(0,0)[lt]{\lineheight{1.25}\smash{\begin{tabular}[t]{l}$k=1$\end{tabular}}}}%
    \put(0,0){\includegraphics[width=\unitlength,page=2]{geodesics.pdf}}%
    \put(0.90922626,0.49102946){\color[rgb]{0,0,0}\makebox(0,0)[lt]{\lineheight{1.25}\smash{\begin{tabular}[t]{l}$k=2$\end{tabular}}}}%
    \put(0,0){\includegraphics[width=\unitlength,page=3]{geodesics.pdf}}%
    \put(0.13631922,0.00726667){\color[rgb]{0,0,0}\makebox(0,0)[lt]{\lineheight{1.25}\smash{\begin{tabular}[t]{l}Original Image\end{tabular}}}}%
    \put(0.62119097,0.00726667){\color[rgb]{0,0,0}\makebox(0,0)[lt]{\lineheight{1.25}\smash{\begin{tabular}[t]{l}Reconstructed Image\end{tabular}}}}%
  \end{picture}%
\endgroup%

%% file: geodesicsTWithVarSpeed.pdf_tex
\begingroup%
  \makeatletter%
  \providecommand\color[2][]{%
    \errmessage{(Inkscape) Color is used for the text in Inkscape, but the package 'color.sty' is not loaded}%
    \renewcommand\color[2][]{}%
  }%
  \providecommand\transparent[1]{%
    \errmessage{(Inkscape) Transparency is used (non-zero) for the text in Inkscape, but the package 'transparent.sty' is not loaded}%
    \renewcommand\transparent[1]{}%
  }%
  \providecommand\rotatebox[2]{#2}%
  \newcommand*\fsize{\dimexpr\f@size pt\relax}%
  \newcommand*\lineheight[1]{\fontsize{\fsize}{#1\fsize}\selectfont}%
  \ifx\svgwidth\undefined%
    \setlength{\unitlength}{360bp}%
    \ifx\svgscale\undefined%
      \relax%
    \else%
      \setlength{\unitlength}{\unitlength * \real{\svgscale}}%
    \fi%
  \else%
    \setlength{\unitlength}{\svgwidth}%
  \fi%
  \global\let\svgwidth\undefined%
  \global\let\svgscale\undefined%
  \makeatother%
  \begin{picture}(1,0.52027802)%
    \lineheight{1}%
    \setlength\tabcolsep{0pt}%
    \put(0,0){\includegraphics[width=\unitlength,page=1]{geodesicsTWithVarSpeed.pdf}}%
    \put(0.13134236,0.00726666){\color[rgb]{0,0,0}\makebox(0,0)[lt]{\lineheight{1.25}\smash{\begin{tabular}[t]{l}Original Image\end{tabular}}}}%
    \put(0.62616785,0.00726666){\color[rgb]{0,0,0}\makebox(0,0)[lt]{\lineheight{1.25}\smash{\begin{tabular}[t]{l}Reconstructed Image\end{tabular}}}}%
    \put(0.30376928,0.2747949){\color[rgb]{0,0,0}\makebox(0,0)[lt]{\lineheight{1.25}\smash{\begin{tabular}[t]{l}$(x,\xi)$\end{tabular}}}}%
    \put(0.7719171,0.26666631){\color[rgb]{0,0,0}\makebox(0,0)[lt]{\lineheight{1.25}\smash{\begin{tabular}[t]{l}$(x,\xi)$\end{tabular}}}}%
    \put(0.54762768,0.26516101){\color[rgb]{0,0,0}\makebox(0,0)[lt]{\lineheight{1.25}\smash{\begin{tabular}[t]{l}$(\tilde{x}_-,\tilde{\xi}_-)$\end{tabular}}}}%
    \put(0.88101227,0.28547379){\color[rgb]{0,0,0}\makebox(0,0)[lt]{\lineheight{1.25}\smash{\begin{tabular}[t]{l}$(\tilde{x}_+,\tilde{\xi}_+)$\end{tabular}}}}%
    \put(0,0){\includegraphics[width=\unitlength,page=2]{geodesicsTWithVarSpeed.pdf}}%
  \end{picture}%
\endgroup%

%% file: bothYV3.pdf_tex
\begingroup%
  \makeatletter%
  \providecommand\color[2][]{%
    \errmessage{(Inkscape) Color is used for the text in Inkscape, but the package 'color.sty' is not loaded}%
    \renewcommand\color[2][]{}%
  }%
  \providecommand\transparent[1]{%
    \errmessage{(Inkscape) Transparency is used (non-zero) for the text in Inkscape, but the package 'transparent.sty' is not loaded}%
    \renewcommand\transparent[1]{}%
  }%
  \providecommand\rotatebox[2]{#2}%
  \newcommand*\fsize{\dimexpr\f@size pt\relax}%
  \newcommand*\lineheight[1]{\fontsize{\fsize}{#1\fsize}\selectfont}%
  \ifx\svgwidth\undefined%
    \setlength{\unitlength}{360bp}%
    \ifx\svgscale\undefined%
      \relax%
    \else%
      \setlength{\unitlength}{\unitlength * \real{\svgscale}}%
    \fi%
  \else%
    \setlength{\unitlength}{\svgwidth}%
  \fi%
  \global\let\svgwidth\undefined%
  \global\let\svgscale\undefined%
  \makeatother%
  \begin{picture}(1,0.52983356)%
    \lineheight{1}%
    \setlength\tabcolsep{0pt}%
    \put(0,0){\includegraphics[width=\unitlength,page=1]{bothYV3.pdf}}%
    \put(0.13192644,0.00726668){\color[rgb]{0,0,0}\makebox(0,0)[lt]{\lineheight{1.25}\smash{\begin{tabular}[t]{l}Original Image\end{tabular}}}}%
    \put(0.6242425,0.00726668){\color[rgb]{0,0,0}\makebox(0,0)[lt]{\lineheight{1.25}\smash{\begin{tabular}[t]{l}Reconstructed Image\end{tabular}}}}%
    \put(0,0){\includegraphics[width=\unitlength,page=2]{bothYV3.pdf}}%
    \put(0.1563488,0.23377284){\color[rgb]{0,0,0}\makebox(0,0)[lt]{\lineheight{1.25}\smash{\begin{tabular}[t]{l}$(x,\xi)$\end{tabular}}}}%
    \put(0.68709036,0.23251265){\color[rgb]{0,0,0}\makebox(0,0)[lt]{\lineheight{1.25}\smash{\begin{tabular}[t]{l}$(x,\xi)$\end{tabular}}}}%
    \put(0.73845085,0.45365585){\color[rgb]{0,0,0}\makebox(0,0)[lt]{\lineheight{1.25}\smash{\begin{tabular}[t]{l}$(\tilde{x},\tilde{\xi})$\end{tabular}}}}%
  \end{picture}%
\endgroup%

%% file: bothVarYV2.pdf_tex
\begingroup%
  \makeatletter%
  \providecommand\color[2][]{%
    \errmessage{(Inkscape) Color is used for the text in Inkscape, but the package 'color.sty' is not loaded}%
    \renewcommand\color[2][]{}%
  }%
  \providecommand\transparent[1]{%
    \errmessage{(Inkscape) Transparency is used (non-zero) for the text in Inkscape, but the package 'transparent.sty' is not loaded}%
    \renewcommand\transparent[1]{}%
  }%
  \providecommand\rotatebox[2]{#2}%
  \newcommand*\fsize{\dimexpr\f@size pt\relax}%
  \newcommand*\lineheight[1]{\fontsize{\fsize}{#1\fsize}\selectfont}%
  \ifx\svgwidth\undefined%
    \setlength{\unitlength}{359.99996948bp}%
    \ifx\svgscale\undefined%
      \relax%
    \else%
      \setlength{\unitlength}{\unitlength * \real{\svgscale}}%
    \fi%
  \else%
    \setlength{\unitlength}{\svgwidth}%
  \fi%
  \global\let\svgwidth\undefined%
  \global\let\svgscale\undefined%
  \makeatother%
  \begin{picture}(1,0.53128946)%
    \lineheight{1}%
    \setlength\tabcolsep{0pt}%
    \put(0,0){\includegraphics[width=\unitlength,page=1]{bothVarYV2.pdf}}%
    \put(0.18366953,0.29337515){\color[rgb]{0,0,0}\makebox(0,0)[lt]{\lineheight{1.25}\smash{\begin{tabular}[t]{l}$(x,\xi)$\end{tabular}}}}%
    \put(0.65245794,0.25147699){\color[rgb]{0,0,0}\makebox(0,0)[lt]{\lineheight{1.25}\smash{\begin{tabular}[t]{l}$(x,\xi)$\end{tabular}}}}%
    \put(0.53564516,0.24027655){\color[rgb]{0,0,0}\makebox(0,0)[lt]{\lineheight{1.25}\smash{\begin{tabular}[t]{l}$(\tilde{x}_+,\tilde{\xi}_+)$\end{tabular}}}}%
    \put(0.74845154,0.24101353){\color[rgb]{0,0,0}\makebox(0,0)[lt]{\lineheight{1.25}\smash{\begin{tabular}[t]{l}$(\tilde{x}_-,\tilde{\xi}_-)$\end{tabular}}}}%
    \put(0.13155463,0.00726667){\color[rgb]{0,0,0}\makebox(0,0)[lt]{\lineheight{1.25}\smash{\begin{tabular}[t]{l}Original Image\end{tabular}}}}%
    \put(0.62542186,0.00726667){\color[rgb]{0,0,0}\makebox(0,0)[lt]{\lineheight{1.25}\smash{\begin{tabular}[t]{l}Reconstructed Image\end{tabular}}}}%
  \end{picture}%
\endgroup%

%% file: wormReconT.pdf_tex
\begingroup%
  \makeatletter%
  \providecommand\color[2][]{%
    \errmessage{(Inkscape) Color is used for the text in Inkscape, but the package 'color.sty' is not loaded}%
    \renewcommand\color[2][]{}%
  }%
  \providecommand\transparent[1]{%
    \errmessage{(Inkscape) Transparency is used (non-zero) for the text in Inkscape, but the package 'transparent.sty' is not loaded}%
    \renewcommand\transparent[1]{}%
  }%
  \providecommand\rotatebox[2]{#2}%
  \newcommand*\fsize{\dimexpr\f@size pt\relax}%
  \newcommand*\lineheight[1]{\fontsize{\fsize}{#1\fsize}\selectfont}%
  \ifx\svgwidth\undefined%
    \setlength{\unitlength}{360bp}%
    \ifx\svgscale\undefined%
      \relax%
    \else%
      \setlength{\unitlength}{\unitlength * \real{\svgscale}}%
    \fi%
  \else%
    \setlength{\unitlength}{\svgwidth}%
  \fi%
  \global\let\svgwidth\undefined%
  \global\let\svgscale\undefined%
  \makeatother%
  \begin{picture}(1,0.5197)%
    \lineheight{1}%
    \setlength\tabcolsep{0pt}%
    \put(0.13168957,0.0072667){\color[rgb]{0,0,0}\makebox(0,0)[lt]{\lineheight{1.25}\smash{\begin{tabular}[t]{l}Original Image\end{tabular}}}}%
    \put(0.637905,0.00916392){\color[rgb]{0,0,0}\makebox(0,0)[lt]{\lineheight{1.25}\smash{\begin{tabular}[t]{l}Reconstructed Image\end{tabular}}}}%
    \put(0,0){\includegraphics[width=\unitlength,page=1]{wormReconT.pdf}}%
  \end{picture}%
\endgroup%

%% file: wormDataT.pdf_tex
\begingroup%
  \makeatletter%
  \providecommand\color[2][]{%
    \errmessage{(Inkscape) Color is used for the text in Inkscape, but the package 'color.sty' is not loaded}%
    \renewcommand\color[2][]{}%
  }%
  \providecommand\transparent[1]{%
    \errmessage{(Inkscape) Transparency is used (non-zero) for the text in Inkscape, but the package 'transparent.sty' is not loaded}%
    \renewcommand\transparent[1]{}%
  }%
  \providecommand\rotatebox[2]{#2}%
  \newcommand*\fsize{\dimexpr\f@size pt\relax}%
  \newcommand*\lineheight[1]{\fontsize{\fsize}{#1\fsize}\selectfont}%
  \ifx\svgwidth\undefined%
    \setlength{\unitlength}{395.02441406bp}%
    \ifx\svgscale\undefined%
      \relax%
    \else%
      \setlength{\unitlength}{\unitlength * \real{\svgscale}}%
    \fi%
  \else%
    \setlength{\unitlength}{\svgwidth}%
  \fi%
  \global\let\svgwidth\undefined%
  \global\let\svgscale\undefined%
  \makeatother%
  \begin{picture}(1,0.54997458)%
    \lineheight{1}%
    \setlength\tabcolsep{0pt}%
    \put(0,0){\includegraphics[width=\unitlength,page=1]{wormDataT.pdf}}%
    \put(0.01700914,0.03839765){\color[rgb]{0,0,0}\makebox(0,0)[lt]{\lineheight{1.25}\smash{\begin{tabular}[t]{l}Collected Data\end{tabular}}}}%
    \put(0.22870508,0.03839765){\color[rgb]{0,0,0}\makebox(0,0)[lt]{\lineheight{1.25}\smash{\begin{tabular}[t]{l}Downsampled Data\end{tabular}}}}%
    \put(0.6481121,0.03839765){\color[rgb]{0,0,0}\makebox(0,0)[t]{\lineheight{1.25}\smash{\begin{tabular}[t]{c}Fourier Transform\\of Collected Data\end{tabular}}}}%
    \put(0.90826093,0.03839765){\color[rgb]{0,0,0}\makebox(0,0)[t]{\lineheight{1.25}\smash{\begin{tabular}[t]{c}Fourier Transform of\\Undersampled Data\end{tabular}}}}%
  \end{picture}%
\endgroup%

%% file: wormReconY.pdf_tex
\begingroup%
  \makeatletter%
  \providecommand\color[2][]{%
    \errmessage{(Inkscape) Color is used for the text in Inkscape, but the package 'color.sty' is not loaded}%
    \renewcommand\color[2][]{}%
  }%
  \providecommand\transparent[1]{%
    \errmessage{(Inkscape) Transparency is used (non-zero) for the text in Inkscape, but the package 'transparent.sty' is not loaded}%
    \renewcommand\transparent[1]{}%
  }%
  \providecommand\rotatebox[2]{#2}%
  \newcommand*\fsize{\dimexpr\f@size pt\relax}%
  \newcommand*\lineheight[1]{\fontsize{\fsize}{#1\fsize}\selectfont}%
  \ifx\svgwidth\undefined%
    \setlength{\unitlength}{360bp}%
    \ifx\svgscale\undefined%
      \relax%
    \else%
      \setlength{\unitlength}{\unitlength * \real{\svgscale}}%
    \fi%
  \else%
    \setlength{\unitlength}{\svgwidth}%
  \fi%
  \global\let\svgwidth\undefined%
  \global\let\svgscale\undefined%
  \makeatother%
  \begin{picture}(1,0.51900549)%
    \lineheight{1}%
    \setlength\tabcolsep{0pt}%
    \put(0.13134235,0.00726665){\color[rgb]{0,0,0}\makebox(0,0)[lt]{\lineheight{1.25}\smash{\begin{tabular}[t]{l}Original Image\end{tabular}}}}%
    \put(0.62616784,0.00726665){\color[rgb]{0,0,0}\makebox(0,0)[lt]{\lineheight{1.25}\smash{\begin{tabular}[t]{l}Reconstructed Image\end{tabular}}}}%
    \put(0,0){\includegraphics[width=\unitlength,page=1]{wormReconY.pdf}}%
  \end{picture}%
\endgroup%

%% file: wormDataY.pdf_tex
\begingroup%
  \makeatletter%
  \providecommand\color[2][]{%
    \errmessage{(Inkscape) Color is used for the text in Inkscape, but the package 'color.sty' is not loaded}%
    \renewcommand\color[2][]{}%
  }%
  \providecommand\transparent[1]{%
    \errmessage{(Inkscape) Transparency is used (non-zero) for the text in Inkscape, but the package 'transparent.sty' is not loaded}%
    \renewcommand\transparent[1]{}%
  }%
  \providecommand\rotatebox[2]{#2}%
  \newcommand*\fsize{\dimexpr\f@size pt\relax}%
  \newcommand*\lineheight[1]{\fontsize{\fsize}{#1\fsize}\selectfont}%
  \ifx\svgwidth\undefined%
    \setlength{\unitlength}{394.5968399bp}%
    \ifx\svgscale\undefined%
      \relax%
    \else%
      \setlength{\unitlength}{\unitlength * \real{\svgscale}}%
    \fi%
  \else%
    \setlength{\unitlength}{\svgwidth}%
  \fi%
  \global\let\svgwidth\undefined%
  \global\let\svgscale\undefined%
  \makeatother%
  \begin{picture}(1,0.55057047)%
    \lineheight{1}%
    \setlength\tabcolsep{0pt}%
    \put(0,0){\includegraphics[width=\unitlength,page=1]{wormDataY.pdf}}%
    \put(0.01685872,0.03843922){\color[rgb]{0,0,0}\makebox(0,0)[lt]{\lineheight{1.25}\smash{\begin{tabular}[t]{l}Collected Data\end{tabular}}}}%
    \put(0.22882183,0.03843922){\color[rgb]{0,0,0}\makebox(0,0)[lt]{\lineheight{1.25}\smash{\begin{tabular}[t]{l}Downsampled Data\end{tabular}}}}%
    \put(0.64767116,0.03843922){\color[rgb]{0,0,0}\makebox(0,0)[t]{\lineheight{1.25}\smash{\begin{tabular}[t]{c}Fourier Transform\\of Collected Data\end{tabular}}}}%
    \put(0.90816195,0.03843922){\color[rgb]{0,0,0}\makebox(0,0)[t]{\lineheight{1.25}\smash{\begin{tabular}[t]{c}Fourier Transform of\\Undersampled Data\end{tabular}}}}%
  \end{picture}%
\endgroup%

%% file: avgTV2.pdf_tex
\begingroup%
  \makeatletter%
  \providecommand\color[2][]{%
    \errmessage{(Inkscape) Color is used for the text in Inkscape, but the package 'color.sty' is not loaded}%
    \renewcommand\color[2][]{}%
  }%
  \providecommand\transparent[1]{%
    \errmessage{(Inkscape) Transparency is used (non-zero) for the text in Inkscape, but the package 'transparent.sty' is not loaded}%
    \renewcommand\transparent[1]{}%
  }%
  \providecommand\rotatebox[2]{#2}%
  \newcommand*\fsize{\dimexpr\f@size pt\relax}%
  \newcommand*\lineheight[1]{\fontsize{\fsize}{#1\fsize}\selectfont}%
  \ifx\svgwidth\undefined%
    \setlength{\unitlength}{446.39998627bp}%
    \ifx\svgscale\undefined%
      \relax%
    \else%
      \setlength{\unitlength}{\unitlength * \real{\svgscale}}%
    \fi%
  \else%
    \setlength{\unitlength}{\svgwidth}%
  \fi%
  \global\let\svgwidth\undefined%
  \global\let\svgscale\undefined%
  \makeatother%
  \begin{picture}(1,0.38475299)%
    \lineheight{1}%
    \setlength\tabcolsep{0pt}%
    \put(0,0){\includegraphics[width=\unitlength,page=1]{avgTV2.pdf}}%
    \put(0.06716881,0.04221606){\color[rgb]{0,0,0}\makebox(0,0)[lt]{\lineheight{1.25}\smash{\begin{tabular}[t]{l}Original Image\end{tabular}}}}%
    \put(0.36113297,0.04221606){\color[rgb]{0,0,0}\makebox(0,0)[lt]{\lineheight{1.25}\smash{\begin{tabular}[t]{l}Wave speed with fast spot\end{tabular}}}}%
    \put(0.85285507,0.03946236){\color[rgb]{0,0,0}\makebox(0,0)[t]{\lineheight{1.25}\smash{\begin{tabular}[t]{c}Reconstructed image with \\data averaged in $t$ variable.\end{tabular}}}}%
  \end{picture}%
\endgroup%

%% file: avgYV2.pdf_tex
\begingroup%
  \makeatletter%
  \providecommand\color[2][]{%
    \errmessage{(Inkscape) Color is used for the text in Inkscape, but the package 'color.sty' is not loaded}%
    \renewcommand\color[2][]{}%
  }%
  \providecommand\transparent[1]{%
    \errmessage{(Inkscape) Transparency is used (non-zero) for the text in Inkscape, but the package 'transparent.sty' is not loaded}%
    \renewcommand\transparent[1]{}%
  }%
  \providecommand\rotatebox[2]{#2}%
  \newcommand*\fsize{\dimexpr\f@size pt\relax}%
  \newcommand*\lineheight[1]{\fontsize{\fsize}{#1\fsize}\selectfont}%
  \ifx\svgwidth\undefined%
    \setlength{\unitlength}{417.60001373bp}%
    \ifx\svgscale\undefined%
      \relax%
    \else%
      \setlength{\unitlength}{\unitlength * \real{\svgscale}}%
    \fi%
  \else%
    \setlength{\unitlength}{\svgwidth}%
  \fi%
  \global\let\svgwidth\undefined%
  \global\let\svgscale\undefined%
  \makeatother%
  \begin{picture}(1,0.40189988)%
    \lineheight{1}%
    \setlength\tabcolsep{0pt}%
    \put(0,0){\includegraphics[width=\unitlength,page=1]{avgYV2.pdf}}%
    \put(0.06754849,0.03632186){\color[rgb]{0,0,0}\makebox(0,0)[lt]{\lineheight{1.25}\smash{\begin{tabular}[t]{l}Original Image\end{tabular}}}}%
    \put(0.34887393,0.03632186){\color[rgb]{0,0,0}\makebox(0,0)[lt]{\lineheight{1.25}\smash{\begin{tabular}[t]{l}Wave Speed with fast spot\end{tabular}}}}%
    \put(0.84543109,0.03505749){\color[rgb]{0,0,0}\makebox(0,0)[t]{\lineheight{1.25}\smash{\begin{tabular}[t]{c}Reconstructed image with \\data averaged in $y$ variable.\end{tabular}}}}%
    \put(0,0){\includegraphics[width=\unitlength,page=2]{avgYV2.pdf}}%
  \end{picture}%
\endgroup%

%% file: antiAliasingZebra.pdf_tex
\begingroup%
  \makeatletter%
  \providecommand\color[2][]{%
    \errmessage{(Inkscape) Color is used for the text in Inkscape, but the package 'color.sty' is not loaded}%
    \renewcommand\color[2][]{}%
  }%
  \providecommand\transparent[1]{%
    \errmessage{(Inkscape) Transparency is used (non-zero) for the text in Inkscape, but the package 'transparent.sty' is not loaded}%
    \renewcommand\transparent[1]{}%
  }%
  \providecommand\rotatebox[2]{#2}%
  \newcommand*\fsize{\dimexpr\f@size pt\relax}%
  \newcommand*\lineheight[1]{\fontsize{\fsize}{#1\fsize}\selectfont}%
  \ifx\svgwidth\undefined%
    \setlength{\unitlength}{414bp}%
    \ifx\svgscale\undefined%
      \relax%
    \else%
      \setlength{\unitlength}{\unitlength * \real{\svgscale}}%
    \fi%
  \else%
    \setlength{\unitlength}{\svgwidth}%
  \fi%
  \global\let\svgwidth\undefined%
  \global\let\svgscale\undefined%
  \makeatother%
  \begin{picture}(1,0.3866957)%
    \lineheight{1}%
    \setlength\tabcolsep{0pt}%
    \put(0,0){\includegraphics[width=\unitlength,page=1]{antiAliasingZebra.pdf}}%
    \put(0.08552718,0.04255073){\color[rgb]{0,0,0}\makebox(0,0)[lt]{\lineheight{1.25}\smash{\begin{tabular}[t]{l}Original Image\end{tabular}}}}%
    \put(0.49988927,0.04255073){\color[rgb]{0,0,0}\makebox(0,0)[t]{\lineheight{1.25}\smash{\begin{tabular}[t]{c}Reconstruction from\\undersampled $y$ data\end{tabular}}}}%
    \put(0.84801454,0.04255073){\color[rgb]{0,0,0}\makebox(0,0)[t]{\lineheight{1.25}\smash{\begin{tabular}[t]{c}Reconstruction using\\anti-aliasing scheme\end{tabular}}}}%
  \end{picture}%
\endgroup%

%% file: zebraTV2.pdf_tex
\begingroup%
  \makeatletter%
  \providecommand\color[2][]{%
    \errmessage{(Inkscape) Color is used for the text in Inkscape, but the package 'color.sty' is not loaded}%
    \renewcommand\color[2][]{}%
  }%
  \providecommand\transparent[1]{%
    \errmessage{(Inkscape) Transparency is used (non-zero) for the text in Inkscape, but the package 'transparent.sty' is not loaded}%
    \renewcommand\transparent[1]{}%
  }%
  \providecommand\rotatebox[2]{#2}%
  \newcommand*\fsize{\dimexpr\f@size pt\relax}%
  \newcommand*\lineheight[1]{\fontsize{\fsize}{#1\fsize}\selectfont}%
  \ifx\svgwidth\undefined%
    \setlength{\unitlength}{360bp}%
    \ifx\svgscale\undefined%
      \relax%
    \else%
      \setlength{\unitlength}{\unitlength * \real{\svgscale}}%
    \fi%
  \else%
    \setlength{\unitlength}{\svgwidth}%
  \fi%
  \global\let\svgwidth\undefined%
  \global\let\svgscale\undefined%
  \makeatother%
  \begin{picture}(1,0.56133337)%
    \lineheight{1}%
    \setlength\tabcolsep{0pt}%
    \put(0,0){\includegraphics[width=\unitlength,page=1]{zebraTV2.pdf}}%
    \put(0.1316896,0.04889998){\color[rgb]{0,0,0}\makebox(0,0)[lt]{\lineheight{1.25}\smash{\begin{tabular}[t]{l}Original Image\end{tabular}}}}%
    \put(0.76668308,0.04889998){\color[rgb]{0,0,0}\makebox(0,0)[t]{\lineheight{1.25}\smash{\begin{tabular}[t]{c}Reconstructed Image\\(Undersampled in $t$)\end{tabular}}}}%
  \end{picture}%
\endgroup%

%% file: zebraDataTV2.pdf_tex
\begingroup%
  \makeatletter%
  \providecommand\color[2][]{%
    \errmessage{(Inkscape) Color is used for the text in Inkscape, but the package 'color.sty' is not loaded}%
    \renewcommand\color[2][]{}%
  }%
  \providecommand\transparent[1]{%
    \errmessage{(Inkscape) Transparency is used (non-zero) for the text in Inkscape, but the package 'transparent.sty' is not loaded}%
    \renewcommand\transparent[1]{}%
  }%
  \providecommand\rotatebox[2]{#2}%
  \newcommand*\fsize{\dimexpr\f@size pt\relax}%
  \newcommand*\lineheight[1]{\fontsize{\fsize}{#1\fsize}\selectfont}%
  \ifx\svgwidth\undefined%
    \setlength{\unitlength}{345.60073471bp}%
    \ifx\svgscale\undefined%
      \relax%
    \else%
      \setlength{\unitlength}{\unitlength * \real{\svgscale}}%
    \fi%
  \else%
    \setlength{\unitlength}{\svgwidth}%
  \fi%
  \global\let\svgwidth\undefined%
  \global\let\svgscale\undefined%
  \makeatother%
  \begin{picture}(1,0.47312159)%
    \lineheight{1}%
    \setlength\tabcolsep{0pt}%
    \put(0,0){\includegraphics[width=\unitlength,page=1]{zebraDataTV2.pdf}}%
    \put(0.10370393,0.04374991){\color[rgb]{0,0,0}\makebox(0,0)[t]{\lineheight{1.25}\smash{\begin{tabular}[t]{c}Collected Data\\(Top edge)\end{tabular}}}}%
    \put(0.33998814,0.04374991){\color[rgb]{0,0,0}\makebox(0,0)[t]{\lineheight{1.25}\smash{\begin{tabular}[t]{c}Undersampled\\Data\end{tabular}}}}%
    \put(0.65975399,0.04222213){\color[rgb]{0,0,0}\makebox(0,0)[t]{\lineheight{1.25}\smash{\begin{tabular}[t]{c}$\mathcal{F}(Mf)$\end{tabular}}}}%
    \put(0.89579827,0.04222213){\color[rgb]{0,0,0}\makebox(0,0)[t]{\lineheight{1.25}\smash{\begin{tabular}[t]{c}$\mathcal{F}(Mf)$\\(undersampled)\end{tabular}}}}%
  \end{picture}%
\endgroup%

%% file: zebraYV2.pdf_tex
\begingroup%
  \makeatletter%
  \providecommand\color[2][]{%
    \errmessage{(Inkscape) Color is used for the text in Inkscape, but the package 'color.sty' is not loaded}%
    \renewcommand\color[2][]{}%
  }%
  \providecommand\transparent[1]{%
    \errmessage{(Inkscape) Transparency is used (non-zero) for the text in Inkscape, but the package 'transparent.sty' is not loaded}%
    \renewcommand\transparent[1]{}%
  }%
  \providecommand\rotatebox[2]{#2}%
  \newcommand*\fsize{\dimexpr\f@size pt\relax}%
  \newcommand*\lineheight[1]{\fontsize{\fsize}{#1\fsize}\selectfont}%
  \ifx\svgwidth\undefined%
    \setlength{\unitlength}{360bp}%
    \ifx\svgscale\undefined%
      \relax%
    \else%
      \setlength{\unitlength}{\unitlength * \real{\svgscale}}%
    \fi%
  \else%
    \setlength{\unitlength}{\svgwidth}%
  \fi%
  \global\let\svgwidth\undefined%
  \global\let\svgscale\undefined%
  \makeatother%
  \begin{picture}(1,0.56067219)%
    \lineheight{1}%
    \setlength\tabcolsep{0pt}%
    \put(0,0){\includegraphics[width=\unitlength,page=1]{zebraYV2.pdf}}%
    \put(0.13134235,0.04893331){\color[rgb]{0,0,0}\makebox(0,0)[lt]{\lineheight{1.25}\smash{\begin{tabular}[t]{l}Original Image\end{tabular}}}}%
    \put(0.76703031,0.04893331){\color[rgb]{0,0,0}\makebox(0,0)[t]{\lineheight{1.25}\smash{\begin{tabular}[t]{c}Reconstructed Image\\(Undersampled in $y$)\end{tabular}}}}%
  \end{picture}%
\endgroup%

%% file: zebraDataYV2.pdf_tex
\begingroup%
  \makeatletter%
  \providecommand\color[2][]{%
    \errmessage{(Inkscape) Color is used for the text in Inkscape, but the package 'color.sty' is not loaded}%
    \renewcommand\color[2][]{}%
  }%
  \providecommand\transparent[1]{%
    \errmessage{(Inkscape) Transparency is used (non-zero) for the text in Inkscape, but the package 'transparent.sty' is not loaded}%
    \renewcommand\transparent[1]{}%
  }%
  \providecommand\rotatebox[2]{#2}%
  \newcommand*\fsize{\dimexpr\f@size pt\relax}%
  \newcommand*\lineheight[1]{\fontsize{\fsize}{#1\fsize}\selectfont}%
  \ifx\svgwidth\undefined%
    \setlength{\unitlength}{345.59949875bp}%
    \ifx\svgscale\undefined%
      \relax%
    \else%
      \setlength{\unitlength}{\unitlength * \real{\svgscale}}%
    \fi%
  \else%
    \setlength{\unitlength}{\svgwidth}%
  \fi%
  \global\let\svgwidth\undefined%
  \global\let\svgscale\undefined%
  \makeatother%
  \begin{picture}(1,0.47312328)%
    \lineheight{1}%
    \setlength\tabcolsep{0pt}%
    \put(0,0){\includegraphics[width=\unitlength,page=1]{zebraDataYV2.pdf}}%
    \put(0.1037043,0.04375007){\color[rgb]{0,0,0}\makebox(0,0)[t]{\lineheight{1.25}\smash{\begin{tabular}[t]{c}Collected Data\\(Top edge)\end{tabular}}}}%
    \put(0.34012203,0.04375007){\color[rgb]{0,0,0}\makebox(0,0)[t]{\lineheight{1.25}\smash{\begin{tabular}[t]{c}Undersampled\\Data\end{tabular}}}}%
    \put(0.65975789,0.04222228){\color[rgb]{0,0,0}\makebox(0,0)[t]{\lineheight{1.25}\smash{\begin{tabular}[t]{c}$\mathcal{F}(Mf)$\end{tabular}}}}%
    \put(0.89579995,0.04222228){\color[rgb]{0,0,0}\makebox(0,0)[t]{\lineheight{1.25}\smash{\begin{tabular}[t]{c}$\mathcal{F}(Mf)$\\(undersampled)\end{tabular}}}}%
  \end{picture}%
\endgroup%